\RequirePackage{fix-cm}
\documentclass[aap]{imsart}
\RequirePackage{amsthm,amsmath,amsfonts,amssymb}
\RequirePackage[numbers]{natbib}
\RequirePackage[colorlinks,citecolor=blue,linkcolor=blue,urlcolor=blue]{hyperref}
\RequirePackage{enumitem,booktabs,array}

\startlocaldefs
\makeatletter\setpkgattr{title}{case}{\@firstofone}\makeatother
\makeatletter\def\journal@name{}\makeatother
\DeclareFontShape{T1}{ptm}{m}{scit}{<->ssub*ptm/m/sc}{}
\theoremstyle{plain}
\newtheorem{theorem}{Theorem}
\newtheorem{lemma}{Lemma}
\newtheorem{corollary}{Corollary}
\newtheorem{proposition}{Proposition}
\theoremstyle{definition}
\newtheorem{definition}{Definition}

\renewcommand{\P}{\mathbb{P}}
\newcommand{\E}{\mathbb{E}}
\newcommand{\TV}{\operatorname{TV}}
\newcommand{\ind}{\mathbf 1}

\renewcommand{\paragraph}[1]{\par\smallskip\noindent\textbf{#1}\ }
\endlocaldefs

\begin{document}

\begin{frontmatter}
\title{Detection and Recovery of Core--Periphery Structure in Random Geometric Graphs}
\runtitle{Core--periphery detection and recovery on the cube}

\begin{aug}
\author[A]{\textnormal{Mingao Yuan}\ead[label=e1]{myuan2@utep.edu}}
\address[A]{\textnormal{Department of Mathematical Sciences, The University of Texas at El Paso, El Paso, Texas 79968, USA, \printead{e1}}}
\end{aug}

\begin{abstract}
Core--periphery (CP) structure, characterized by a dense core of interconnected nodes linked to sparsely connected peripheral nodes, is a prevalent feature of many real-world networks. We propose a core--periphery random geometric graph (CP-RGG) model in which vertices are independently and uniformly distributed on the unit cube \([0,1]^m\), and each vertex is designated as a core node with probability \(p_n\). Two vertices are connected if their Euclidean distance falls below a type-dependent threshold, governed by three distinct radii for core--core, core--periphery, and periphery--periphery pairs. We investigate the information-theoretic limits for both the detection and recovery of the core in the CP-RGG. For detection in the non-degenerate regime, we prove that the total variation distance between the CP-RGG and a null RGG without CP structure converges to one if and only if the expected number of core nodes tends to infinity. For recovery, we demonstrate that almost exact recovery is achievable if and only if the expected degree tends to infinity. Finally, we establish limits for exact recovery. In the one-dimensional sparse regime (\(m=1\)), we obtain a sharp threshold when the core--core and core--periphery radii coincide. For \(m\ge2\), when the core has polynomial size and the core--periphery radius exceeds twice the periphery--periphery radius, we prove that the exact recovery threshold is sharp up to constant factors, scaling on the order of \(\log n\).
\end{abstract}

\begin{keyword}[class=MSC]
\kwd{60K35}
\kwd{05C80}
\kwd{62G10}
\end{keyword}

\begin{keyword}
\kwd{random geometric graph}
\kwd{core--periphery structure}
\kwd{total variation distance}
\kwd{induced star}
\kwd{detection threshold}
\end{keyword}

\end{frontmatter}

\section{Introduction}\label{sec:intro}

Core--periphery (CP) structure, in which a small set of well-connected core vertices is surrounded by a periphery that attaches mostly to the core, is a common feature of real networks \cite{BE00,CLWU13,RPFM14,ZMN15,CRLP16,YS23}. In global trade networks, large economies trade with economies of all sizes and form the core, while small economies trade little with each other and form the periphery \cite{Kru96,MLS16}. Airport networks show the same pattern, with major hubs linked to each other and to regional airports \cite{LS17,LS19}, and so do citation networks, in which high-profile papers are cited widely while obscure papers rarely cite each other \cite{SCL20,WPKWC22}. Studying the CP structure of a network is important because core nodes are likely to be more influential: the failure of a core power plant, for instance, has far greater consequences for a power grid than the failure of a peripheral one \cite{YS23}.

A core--periphery structure is fundamentally distinct from the widely studied community structure \cite{ACV14, MNS15, Abb18}, and this distinction shapes both modeling choices and analytical methods \cite{ZMN15,CRLP16,YS23,YSM25}. In a community structure, vertices are partitioned into groups such that edges occur more frequently \emph{within} a group than \emph{between} groups \cite{ACV14, MNS15, Abb18}. To make this precise, consider a two-block stochastic block model \cite{Abb18} with connection probabilities \(q_{11}\), \(q_{10}\), and \(q_{00}\), where \(q_{11}\) denotes the probability of an edge between two vertices in block~1, \(q_{00}\) denotes the probability between two vertices in block~0, and \(q_{10}\) denotes the probability between a vertex in block~1 and a vertex in block~0. An assortative community structure satisfies \(q_{11}, q_{00} > q_{10}\), ensuring that each group is internally dense. By contrast, a core--periphery structure obeys the ordering \(q_{11}\ge q_{10}\ge q_{00}\): block~1 forms the dense core, whereas block~0 forms the \emph{sparsest} region of the graph, meaning that peripheral vertices connect predominantly to the core rather than to one another \cite{ZMN15,CRLP16,YS23,YSM25}. The two blocks play asymmetric roles: the core acts as a set of hubs rather than a conventional community, while the periphery does not constitute a community at all, as its members are more strongly connected to the core than to each other. In addition, standard community detection algorithms, such as modularity maximisation or spectral clustering, are designed to identify internally dense clusters rather than a core. Consequently, the analyses in \cite{ZMN15,CRLP16,YS23,YSM25} are formulated specifically for CP models and estimators.

A substantial body of literature addresses the recovery of core--periphery structure. Zhang, Martin, and Newman \cite{ZMN15} fit a block model with CP-shaped connection probabilities via maximum likelihood using the expectation--maximisation (EM) algorithm. Cucuringu, Rombach, Lee, and Porter \cite{CRLP16} propose spectral and geodesic-path methods to partition vertices into a core and a periphery, while Kojaku and Masuda \cite{KM17} propose a scalable algorithm that detects multiple non-overlapping core--periphery pairs and evaluates the statistical significance of the detected CP structures against randomised null networks. Miao and Li \cite{ML23} investigate the identification of an \emph{informative core}---under the assumption that the periphery carries no information about the core---and demonstrate that a spectral estimator achieves consistent recovery. More recently, Yanchenko and Sengupta \cite{YS26} introduce a rapid label-switching algorithm to optimise the Borgatti--Everett objective, rendering the ideal-structure formulation scalable to large networks. From a statistical inference perspective, Yanchenko, Sengupta, and Mukherjee \cite{YSM25} define a population-level CP-strength parameter and formulate four random-graph models with a planted core, where edge probabilities decrease monotonically across core--core, core--periphery, and periphery--periphery pairs. They establish consistency guarantees for exact core recovery and construct hypothesis tests for the presence of CP structure; additional surveys and developments appear in \cite{CLWU13,YS23}. Crucially, in all of the statistical models above, edges are conditionally independent given the latent node labels, much like in an Erd\H{o}s--R\'enyi graph with a planted dense subgraph \cite{ACV14} or a standard stochastic block model \cite{Abb18}.

Many real-world networks, however, are \emph{spatial}: vertices occupy positions in an underlying metric space, and edges form preferentially between nearby nodes. A canonical model for such networks is the random geometric graph (RGG), in which vertices are realised as random points in space and two vertices are connected whenever their distance falls below a prescribed threshold. Because edges are determined by the same underlying vertex positions, they are not independent; in particular, the neighbours of a given vertex are more likely to be connected to one another than in graphs with independent edges. The RGG originates with Gilbert \cite{Gil61}, who introduced it under the name of random plane networks to model wireless communication, and its foundational theory is systematically developed in the monograph by Penrose \cite{Penrose03}. Since their introduction, RGGs have found widespread applications across diverse disciplines, including statistical physics \cite{DC02} and computational biology \cite{HRP08}. An extensive body of literature continues to explore various facets of random geometric graphs \cite{GRK05,Zha15,BDER16,Yuan24,YY26,Yuan26}; for a comprehensive overview of recent developments, we refer the reader to \cite{DDC23}.

To our knowledge, the detection and recovery of core--periphery structure in random geometric graphs have not yet been investigated. In this paper, we address this gap by introducing the core--periphery random geometric graph (CP-RGG) model and establishing the information-theoretic thresholds---sharp, or sharp up to constant factors---for both tasks. Specifically, \(n\) nodes are distributed independently and uniformly on the unit cube \([0,1]^m\), and each node is independently designated as a core node with probability \(p_n\), with the remaining nodes forming the periphery. Two nodes are connected if their Euclidean distance falls below a type-dependent threshold: core--core, core--periphery, and periphery--periphery pairs connect within radii \(r_np_{11}\), \(r_np_{10}\), and \(r_np_{00}\), respectively (Definition~\ref{def:model}). Here, \(1>p_{11}\ge p_{10}>p_{00}>0\) are fixed constants, while \(r_n\in(0,1)\) may vary with \(n\). The CP-RGG serves as a spatial analogue of the CP stochastic block model of \cite{ZMN15,CRLP16,YS23,YSM25}.

We investigate two fundamental questions for the CP-RGG. The first is \emph{detection}, which we formulate as a hypothesis testing problem following \cite{ACV14,BDER16,MNS15,Abb18}. Under the null hypothesis, the graph is a homogeneous RGG governed by a common connection radius, calibrated so that its expected edge count matches that of the CP-RGG; this null model exhibits no CP structure. Under the alternative hypothesis, the graph is generated from the CP-RGG. Because the expected edge densities are matched, the total edge count carries no statistical signal, ensuring that any discrepancy between the two distributions arises solely from the underlying CP structure. Following Bubeck, Ding, Eldan, and R\'acz \cite{BDER16}, we quantify this discrepancy via the total variation distance. Recall that the total variation distance converges to one if and only if there exists a test whose type~I and type~II error probabilities both vanish asymptotically. We establish the sharp threshold at which this convergence occurs, saying that detection is \emph{possible} if the total variation distance tends to one, and \emph{impossible} if it remains bounded away from one. The second question is \emph{recovery}: observing only the unlabelled graph---without access to vertex positions or latent labels---we determine when the core can be reliably estimated. An estimator achieves \emph{exact recovery} if it coincides with the true core with probability tending to one, and \emph{almost exact recovery} if it coincides with the true core up to a vanishing fraction of the core size with probability tending to one. We say that exact (respectively, almost exact) recovery is \emph{possible} if it is attained by some estimator, and \emph{impossible} otherwise. Both criteria are formulated rigorously in Section~\ref{sec:perf}.

We answer both questions for the CP-RGG, and the results are summarised in Tables~\ref{tab:summary} and~\ref{tab:recovery}. For detection, we show that under mild conditions (the standing assumption~\eqref{eq:B} below), the total variation distance tends to one if and only if the expected number of core nodes, \(np_n\), tends to infinity (Theorem~\ref{thm:main}; Table~\ref{tab:summary}). For almost exact recovery, we prove that under the same conditions, recovery is possible if and only if the expected degree of the CP-RGG, which scales on the order of \(nr_n^m\), tends to infinity (Theorems~\ref{thm:recovery} and~\ref{thm:sharp}; Table~\ref{tab:recovery}). We also establish thresholds for exact recovery, which are sharp or sharp up to constant factors across broad parameter regimes. For \(m=1\), \(p_{11}=p_{10}\), and \(np_nr_n\to0\), the threshold is sharp: for any fixed \(\varepsilon\in(0,1)\), exact recovery is achievable if
\[
nr_n\ge(1+\varepsilon)\frac{\log(np_n)}{2(p_{10}-p_{00})}
\]
and impossible if
\[
nr_n\le(1-\varepsilon)\frac{\log(np_n)}{2(p_{10}-p_{00})} .
\]
For general \(m\), when the core has polynomial size, the limits for exact recovery are sharp up to constant factors. Specifically, if \(p_n=n^{-t}\) for a fixed constant \(t\in(0,1)\) and \(p_{10}>2p_{00}\), then exact recovery is achievable if \(nr_n^m\ge K_1\log n\) and impossible if \(nr_n^m\le K_2\log n\), for some constants \(K_1>K_2>0\). Because the core size in real-world networks is typically much smaller than the total network size \(n\), assuming a polynomial core size is a natural choice that covers a wide range of CP structures \cite{YSM25}.

\begin{table}[t]
\centering
\caption{Behavior of \(\TV(H_0,H_1)\). Here, $K$ is a positive constant.}\label{tab:summary}
\small
\begin{tabular}{@{}p{3cm}p{5cm}p{2.5cm}@{}}
\toprule
Regime & Total variation distance & Detection\\
\midrule
\(np_n\to\infty\) & \(\TV(H_0,H_1)=1+o(1)\) & Possible  
\\
\(np_n=O(1)\) & \(\TV(H_0,H_1)\le1-e^{-K}\) & Impossible\\
\bottomrule
\end{tabular}
\medskip

\centering
\caption{Recovery as \(r_n\to0\). Here, \(t,\varepsilon\in(0,1)\), and \(K_1>K_2>0\) are fixed constants.
}\label{tab:recovery}
\small
\renewcommand{\arraystretch}{1}
\begin{tabular}{@{}>{\raggedright\arraybackslash}p{8.3cm}>{\raggedright\arraybackslash}p{5.4cm}@{}}
\toprule
Regime & Recovery\\
\midrule
\(nr_n^m\to\infty\) & Almost exact recovery is possible\\
\(nr_n^m=O(1)\) & Almost exact recovery is impossible\\
\midrule

$m=1$, $p_{11}=p_{10}$, \(np_nr_n\to0\), $nr_n\ge(1+\varepsilon)\frac{\log(np_n)}{2(p_{10}-p_{00})}$ & Exact recovery is possible\\ [8pt]
$m=1$, $p_{11}=p_{10}$, \(np_nr_n\to0\), $nr_n\le(1-\varepsilon)\frac{\log(np_n)}{2(p_{10}-p_{00})}$ & Exact recovery is impossible\\ [8pt]
\midrule
\(p_n=n^{-t}\), \(p_{10}>2p_{00}\),  
\(nr_n^m\ge K_1\log n\) &
Exact recovery is possible\\[8pt]
\(p_n=n^{-t}\), \(p_{10}>2p_{00}\),  
\(nr_n^m\le K_2\log n\) &
Exact recovery is impossible \\
\bottomrule
\end{tabular}
\end{table}

The remainder of the paper is organised as follows. Section~\ref{sec:model} fixes the notation, defines the graph models and compares them with related models, and introduces the hypotheses, the total variation distance and the notions of recovery. Section~\ref{sec:results} states the detection and recovery results, which are summarised in Tables~\ref{tab:summary} and~\ref{tab:recovery}. Section~\ref{sec:proofs} contains the proofs.

\section{Models and hypotheses}\label{sec:model}

This section fixes the models. We first recall the notation (Section~\ref{sec:notation}), then define the two graph models (Section~\ref{sec:graphs}), the hypotheses and the total variation distance used to measure detection (Section~\ref{sec:hyp}), and the notions of recovery (Section~\ref{sec:perf}).

\subsection{Notation}\label{sec:notation}

Throughout, \(\|x\|\) denotes the Euclidean norm of \(x\in\mathbb R^m\),
\[
V_m=\frac{\pi^{m/2}}{\Gamma(m/2+1)}
\]
is the volume of the unit ball, and, for \(0<\rho\le1\), \(f(\rho)\) is the \emph{edge-probability function}
\[
f(\rho)=\P\bigl(\|X-Y\|<\rho\bigr),\qquad X,Y\ \text{independent and uniform on }[0,1]^m .
\]
The function \(f\) has no closed form; Lemma~\ref{lem:edgeprob} gives the bounds that are needed. We also need a geometric constant: let \(s_m\) be the largest number of unit vectors in \(\mathbb R^m\) whose pairwise angles all exceed \(\pi/3\). It is finite, it is at most the kissing number of \(\mathbb R^m\), and \(s_1=2\), \(s_2=5\) (the regular pentagon; six directions would require an angle of \(\pi/3\) between neighbours) and \(s_3=12\) (the icosahedron, whose angles are about \(63.4^\circ\)).

We also use the Bachmann--Landau notation. For positive sequences \((a_n)\) and \((b_n)\), we write \(a_n=O(b_n)\) if \(\limsup_{n\to\infty}a_n/b_n<\infty\), and \(a_n=o(b_n)\) if \(a_n/b_n\to0\). We write \(a_n\ll b_n\) if \(a_n=o(b_n)\), \(a_n\gg b_n\) if \(b_n=o(a_n)\), and \(a_n\gtrsim b_n\) if \(a_n\ge cb_n\) for some constant \(c>0\) and all large \(n\). We say that \(a_n\) is \emph{of order} \(b_n\) if both \(a_n=O(b_n)\) and \(b_n=O(a_n)\). An event holds \emph{with high probability} if its probability tends to one as \(n\to\infty\).

\subsection{The two graph models}\label{sec:graphs}

\begin{definition}[Homogeneous RGG on the cube]\label{def:hom}
Let \(m\ge1\) and \(n\ge2\) be integers and let \(\rho\in(0,1)\). Let \(X_1,\dots,X_n\) be i.i.d.\ uniform on the cube \([0,1]^m\), with no identification of opposite faces. The \emph{homogeneous random geometric graph} \(\mathcal G(\rho)\) is the undirected graph on the vertex set \(\{1,\dots,n\}\) with adjacency matrix \(A\) given by \(A_{ii}=0\) and
\[
A_{ij}=\ind\{\|X_i-X_j\|<\rho\},\qquad 1\le i\ne j\le n .
\]
\end{definition}

The graph \(\mathcal G(\rho)\) is the random geometric graph on the cube \([0,1]^m\) with the Euclidean distance and without wrap-around. It is a version of the random geometric graph introduced in \cite{Gil61,Penrose03} and a standard model in both the theory and the applications of random geometric graphs. For example, on the theoretical side, Goel, Rai and Krishnamachari~\cite{GRK05} show that monotone graph properties have sharp thresholds in this model, and Zhao~\cite{Zha15} studies the threshold for the absence of isolated nodes when \(m=1\). On the applied side, the model has been fitted to protein--protein interaction networks~\cite{HRP08} and used in sampling-based motion planning~\cite{SSH18}.

The graph has the following features. First, all pairs of vertices in \(\mathcal G(\rho)\) use the same connection radius \(\rho\), the edges are determined by the positions only, and the edge probability of a pair is \(f(\rho)\). Second, the degree depends on the position: a vertex at distance at least \(\rho\) from the boundary has the whole ball inside the cube, so that each other vertex falls in it with probability \(V_m\rho^m\), while a vertex near the boundary has fewer neighbours. The expected degree of a vertex is therefore of order \(n\rho^m\). Third, the edges are not independent: if \(i\sim j\) and \(j\sim k\), then \(i\) and \(k\) are within distance \(2\rho\), so that the neighbours of a vertex tend to be linked to each other. Finally, the boundary of the cube is not wrapped around, so boundary effects are present.

Building on the random geometric graph \(\mathcal G(\rho)\), we propose a core--periphery random geometric graph. The model retains the vertex positions, the Euclidean norm and the cube of \(\mathcal G(\rho)\), and modifies only the connection radius of a pair, which now depends on the labels of its two vertices.

\begin{definition}[Core--periphery RGG on the cube]\label{def:model}
Let \(m\ge1\) and \(n\ge2\) be integers, \(1>p_{11}\ge p_{10}>p_{00}>0\) be constants, and \(r_n\in(0,1)\) and \(p_n\in(0,1)\) be parameters. Let \(Z_1,\dots,Z_n\) be i.i.d.\ Bernoulli\((p_n)\) random variables and \(X_1,\dots,X_n\) be i.i.d.\ uniform random vectors on \([0,1]^m\), independent of \(Z\). The \emph{core--periphery random geometric graph} (CP-RGG) \(\mathcal G(r_n,p_n,p_{11},p_{10},p_{00})\) is the undirected graph on \(\{1,\dots,n\}\) with adjacency matrix given by \(A_{ii}=0\) and, for \(i\ne j\),
\[
A_{ij}=\ind\{\|X_i-X_j\|<r_np_{Z_iZ_j}\},
\]
where \(p_{Z_iZ_j}\) equals \(p_{11}\), \(p_{10}\) or \(p_{00}\) according to whether \(Z_i+Z_j=2\), \(1\) or \(0\).
\end{definition}

In the CP-RGG \(\mathcal G(r_n,p_n,p_{11},p_{10},p_{00})\), a vertex with \(Z_i=1\) is called a \emph{core} vertex, and a vertex with \(Z_i=0\) a \emph{peripheral} vertex. We treat these vertex labels as random, analogous to the community-structure setting \cite{Abb18,MNS15}. Since the labels are independent of the positions, the core is spread uniformly over the cube, and its size \(|C|\) follows a Binomial\((n,p_n)\) distribution with mean \(np_n\); this is the number of core vertices that appears in the thresholds below. The constants \(p_{11}\), \(p_{10}\), and \(p_{00}\) are the connection ranges of a pair relative to the base radius \(r_n\), which controls the sparsity of the graph: a larger \(r_n\) produces more edges. The condition \(p_{10}>p_{00}\) is what makes the labels visible. If \(p_{10}=p_{00}\), a core vertex and a peripheral vertex have the same range towards peripheral vertices, and the only effect of the core is a longer range between core vertices. The case \(p_{11}=p_{10}=p_{00}\) is the homogeneous graph, in which the core has no effect. When \(p_{11}=p_{10}\), the CP-RGG provides a random model for a noisy version of the ``ideal'' core--periphery structure, in which core vertices are connected to all other vertices and peripheral vertices are mutually disconnected \cite{BE00,YSM25}. Finally, the bounds \(p_{11}<1\) and \(r_n<1\) keep every radius below \(1\), the side of the cube, so that all the edge probabilities are given by the function \(f\) of Section~\ref{sec:notation}. We write the connection radius as the product \(r_np_{Z_iZ_j}\), where \(r_n\) controls the sparsity of the network and \(p_{Z_iZ_j}\) encodes the CP structure, as is customary in the existing literature (see, for example, \cite{YSM25}).

Definition~\ref{def:model} is the geometric analogue of the core--periphery stochastic block model studied in \cite{ZMN15,CRLP16,YS23,YSM25}. In a two-block stochastic block model, each vertex carries a latent label, and the probability of an edge between two vertices depends only on their labels, through connection probabilities \(q_{11}\), \(q_{10}\) and \(q_{00}\): \(q_{11}\) is the probability of an edge between two vertices in block~1, \(q_{00}\) is the probability of an edge between two vertices in block~0, and \(q_{10}\) is the probability of an edge between the two blocks. A CP stochastic block model in the sense of \cite{ZMN15,CRLP16,YSM25} has \(q_{11}\ge q_{10}> q_{00}\) and, in the strict Borgatti--Everett form \cite{BE00}, \(q_{11}=1\) and \(q_{00}=0\). In Definition~\ref{def:model}, the labels \(Z_i\) play the same role, but they determine the connection radius of a pair, \(r_np_{11}\), \(r_np_{10}\) or \(r_np_{00}\), instead of a connection probability. The marginal probability that two vertices are linked is \(f(r_np_{Z_iZ_j})\), which is increasing in the connection radius. Consequently, the ordering \(p_{11}\ge p_{10}>p_{00}\) of the radii corresponds to the ordering of the three connection probabilities required of a CP structure \cite{CRLP16,YSM25}. The difference is that, given the labels, the edges are not independent: they are driven by the common positions \(X_i\) rather than by independent coin tosses.

Definition~\ref{def:model} also differs from the geometric block model of Galhotra, Mazumdar, Pal, and Saha~\cite{GMPS23}. In that model, the vertices have latent community labels and positions, and two vertices are linked when their distance is below a radius that depends on whether their labels agree: the radius is longer within a community than between communities. For this model, they study the limits for recovering the communities. Our model differs in two ways. First, in the geometric block model the radius of a pair depends on whether the two labels agree, so that the labels are exchangeable and each group is dense inside; here the radius depends on \(Z_i+Z_j\), so that the two groups play different roles: the core reaches farthest and the periphery is the sparsest part. Second, the goal there is to recover communities, each of which contains half of the vertices, whereas here we also ask whether a core is present at all, and the core may be a vanishing fraction of the vertices; detection and recovery are therefore separate problems with different thresholds.

Definition~\ref{def:model} also differs from the botnet model of Bet et al.~\cite{BBCH21}. In the botnet model, the vertices outside the botnet form a homogeneous random geometric graph on the torus, while each botnet vertex ignores the positions and connects to every other vertex independently, with a fixed probability equal to the edge probability of a pair of non-botnet vertices. In our model, by contrast, no vertex ignores the geometry: every vertex, core or peripheral, is linked to the vertices within its radius, and a core vertex differs from a peripheral vertex only through a longer radius.

\subsection{Hypotheses}\label{sec:hyp}

As in the detection problems for community structure or geometry in networks that have been widely studied (see, for example, \cite{ACV14,MNS15,BDER16,Abb18,YY26}), we formulate the detection of a CP structure as a hypothesis testing problem.

A pair of vertices contains two core vertices with probability \(p_n^2\), exactly one core vertex with probability \(2p_n(1-p_n)\), and no core vertex with probability \((1-p_n)^2\). Hence the edge probability of \(\mathcal G(r_n,p_n,p_{11},p_{10},p_{00})\) is
\[
p_n^2f(r_np_{11})+2p_n(1-p_n)f(r_np_{10})+(1-p_n)^2f(r_np_{00}).
\]
The \emph{density-matched radius} \(\gamma_n\) is the unique solution in \([r_np_{00},r_np_{11}]\) of
\[
f(\gamma_n)=p_n^2f(r_np_{11})+2p_n(1-p_n)f(r_np_{10})+(1-p_n)^2f(r_np_{00}).
\]
It exists and is unique because \(f\) is continuous and strictly increasing on \([0,1]\) (Lemma~\ref{lem:edgeprob}).

We consider the hypotheses
\[
H_0:\ A\sim\mathcal G(\gamma_n)\qquad\text{against}\qquad H_1:\ A\sim\mathcal G(r_n,p_n,p_{11},p_{10},p_{00}).
\]
Under \(H_0\), the graph is a homogeneous random geometric graph without a core--periphery structure. Under \(H_1\), the graph is a CP-RGG and hence has a core--periphery structure.

The choice of \(\gamma_n\) ensures that \(H_0\) and \(H_1\) have the same expected edge density. Without this matching, a test based solely on the edge count would detect the core whenever \(p_n\) is not too small, simply because the core adds edges. The density-matched null removes this uninformative signal, so that any evidence for the core must come from the CP structure. The same matching is used in block-model detection, where the null has the same expected number of edges as the planted model \cite{MNS15,ACV14,Abb18}, and in the botnet model, where the botnet vertices receive the same expected degree as the other vertices \cite{BBCH21}. The matching is exact for the first moment only: the alternative is a mixture of three radii, and this mixture leaves information in the local structure of the graph.

Detection is a testing problem between the null and the alternative, and how well any test can perform is governed by the total variation distance between the two laws of the graph.

\begin{definition}[Total variation distance]\label{def:tv}
Let \(P\) and \(Q\) be probability measures on a measurable space \((\Omega,\mathcal F)\). Their \emph{total variation distance} is
\[
\TV(P,Q)=\sup_{B\in\mathcal F}\bigl|P(B)-Q(B)\bigr|\in[0,1].
\]
\end{definition}

If \(P\) and \(Q\) have densities \(p\) and \(q\) with respect to a common dominating \(\sigma\)-finite measure \(\mu\) (for instance, \(\mu=P+Q\)), then
\[
\TV(P,Q)=\frac12\int_\Omega|p-q|d\mu=1-\int_\Omega\min(p,q)d\mu .
\]
In particular, when \(\Omega\) is finite or countable and \(\mu\) is the counting measure,
\[
\TV(P,Q)=1-\sum_{\omega\in\Omega}\min\{P(\omega),Q(\omega)\} .
\]

The total variation distance measures how well two models can be told apart: for every measurable test \(\psi:\Omega=\{0,1\}^{\binom n2}\to\{0,1\}\) between \(H_0:A\sim P\) and \(H_1:A\sim Q\),
\[
P(\psi=1)+Q(\psi=0)\ \ge\ 1-\TV(P,Q),
\]
Hence there exists a test whose type~I and type~II error probabilities both tend to zero if and only if \(\TV(P_n,Q_n)\to1\), and no such test exists if \(\limsup_{n\to\infty}\TV(P_n,Q_n)<1\).

\subsection{Recovery of the core}\label{sec:perf}

Recovery is an estimation problem: from the graph alone we must name the core, and we measure the quality of an estimate by the number of vertices it misclassifies, compared with the size of the core.

An estimator of the core is a map \(\widehat C\) from adjacency matrices to subsets of \(\{1,\dots,n\}\); it may depend on \(A\) and on the known parameters, but not on the positions \(X\) or the labels \(Z\), which are not observed. Let \(C=\{i:Z_i=1\}\) be the core and let \(\triangle\) denote the symmetric difference. The notions of recovery follow those used for community detection in block models \cite{ABH16,Abb18,MNS15} and for core--periphery models in the work of Yanchenko, Sengupta, and Mukherjee \cite{YSM25}.

\begin{definition}[Exact and almost exact recovery]\label{def:recovery}
An estimator \(\widehat C=\widehat C_n\) achieves
\begin{enumerate}[label=\textup{(\alph*)},itemsep=1pt]
\item \emph{exact recovery} if, as \(n\to\infty\),
\[
\P(\widehat C=C)\to1 ;
\]
\item \emph{almost exact recovery} if \(|\widehat C\triangle C|/|C|\to0\) in probability, that is, if for every \(\varepsilon>0\),
\[
\P\bigl(|\widehat C\triangle C|\ge\varepsilon|C|\bigr)\to0 .
\]
\end{enumerate}
\end{definition}

Exact (respectively, almost exact) recovery is said to be \emph{possible} if there exists an estimator that achieves it, and \emph{impossible} if no estimator achieves it. In part~(b), the number of misclassified vertices is compared with the size \(|C|\) of the core rather than with \(n\): when \(p_n\to0\), an error of order \(np_n\) is negligible relative to \(n\) but comparable to the size of the entire core. Clearly, exact recovery implies almost exact recovery, so exact recovery generally requires stronger conditions.

\section{Results}\label{sec:results}

Throughout the paper we work under the following standing assumption:
\begin{equation}\label{eq:B}
\liminf_{n\to\infty}nr_n^m>0,\qquad np_n\ge1\ \text{ for all }n,\qquad r_np_{11}\le\tfrac14,\qquad p_n\le\tfrac12 .
\end{equation}
The first two conditions remove degenerate regimes. If \(nr_n^m\to0\), the average degree vanishes, most vertices are isolated, and the labels leave almost no trace. If \(np_n\to0\), then with probability at least \((1-p_n)^n\ge e^{-2np_n}\to1\) there is no core vertex. The constant \(1\) in \(np_n\ge1\) is only a normalisation: any fixed positive lower bound would do. The last two conditions are geometric. The first, \(r_np_{11}\le\tfrac14\), keeps the balls of the largest radius small relative to the cube, and the second, \(p_n\le\tfrac12\), ensures that the core is the smaller of the two groups. Both restrictions are natural, since real networks are typically sparse and the core is typically small \cite{YSM25}.

\subsection{Detection of the core}\label{sec:detection}

The first result is a sharp threshold for detection against the matched null. It has two parts: the two laws are perfectly separated when the expected number of core vertices diverges, and they overlap when it is bounded.

\begin{theorem}\label{thm:main}
Assume~\eqref{eq:B}, and let \(\TV(H_0,H_1)\) denote the total variation distance between the laws of the graph under \(H_0\) and under \(H_1\).
\begin{enumerate}[label=\textup{(\roman*)},itemsep=2pt]
\item \textup{(Detection possible.)} If \(np_n\to\infty\), then
\[
\TV(H_0,H_1)\to1 .
\]
\item \textup{(Detection impossible.)} If \(np_n\le K<\infty\), then
\[
\TV(H_0,H_1)\ \le\ 1-e^{-C_*K},
\]
where
\[
C_*=2+2\left(\frac43\right)^m\left(\left(\frac{p_{11}}{p_{00}}\right)^m-1\right)>0 .
\]
\end{enumerate}
\end{theorem}

Theorem~\ref{thm:main} gives a sharp threshold. Under the standing assumption~\eqref{eq:B}, parts~(i) and~(ii) together say that
\[
\TV(H_0,H_1)\to1\quad\Longleftrightarrow\quad np_n\to\infty ,
\]
so that, once the expected degree is bounded away from zero, the number of core vertices is all that detection requires.

Part~(i) follows from the test that rejects \(H_0\) if and only if the graph contains an induced star \(K_{1,s_m+1}\), that is, a vertex with \(s_m+1\) pairwise non-adjacent neighbours. This test has no type~I error, since a homogeneous random geometric graph never contains such a star, and its type~II error tends to \(0\) when \(np_n\to\infty\). Even when the expected degree is bounded, the test requires only that one of the roughly \(np_n\) core vertices have \(s_m+1\) pairwise separated peripheral neighbours, an event whose probability is bounded below. Thus the threshold is attained even though the matched null removes the information carried by the edge count: what remains is the existence of a core vertex whose peripheral neighbours are spread over a window wider than the periphery radius, which is precisely what the star test detects.

The lower bound in part~(ii) has a different source. If \(np_n\) is bounded, then with probability bounded away from \(0\) the alternative has no core vertex, and it is then a homogeneous graph with the periphery radius. This radius is at relative distance \(O(p_n)=O(1/n)\) from the null radius, and scaling shows that so small a change of radius cannot be detected.

Finally, the constants enter only mildly. Only \(p_{10}>p_{00}\) is used by the star test, through the gap \(p_{10}-p_{00}\) that leaves room for \(s_m+1\) pairwise non-adjacent peripheral neighbours of a core vertex. The constant \(p_{11}\) enters only through the constant \(C_*\) of part~(ii), so \(p_{11}=p_{10}\) is allowed. The third bound in~\eqref{eq:B} serves only to keep all radii at most \(\frac14\), and the proofs go through for any radii at most \(\frac14\). Part~(ii) needs no condition on \(nr_n^m\).

\subsection{Recovery of the core}

Let \(C\) denote the set of core vertices, that is,
\[
C=\{i\in\{1,2,\dots,n\}: Z_i=1\}.
\]
Let \(N=|C|\) and let \(\widehat C_{\mathrm{star}}\) be the set of centres of induced stars \(K_{1,s_m+1}\), that is,
\[
\widehat C_{\mathrm{star}}=\{i\in\{1,2,\dots,n\}:\ i\ \text{has \(s_m+1\) pairwise non-adjacent neighbours}\}.
\]
The set \(\widehat C_{\mathrm{star}}\) can be computed from the graph alone, and the star test of Section~\ref{sec:detection} rejects if and only if \(\widehat C_{\mathrm{star}}\ne\emptyset\). Let \(I=[r_np_{10},1-r_np_{10}]^m\) be the set of positions at distance at least \(r_np_{10}\) from every face of the cube; note that \(r_np_{10}<\frac14\).

The stars find the core vertices in \(I\). A core vertex within distance \(r_np_{10}\) of a face has a truncated neighbourhood, and the \(s_m+1\) balls of the packing construction may not fit. Such a vertex is, however, still better connected than a typical peripheral vertex when the part of its ball inside the cube is large enough, and the simplest statistic that detects this is the degree. For \(\eta\in(0,1)\), let \(\deg(i)\) be the degree of vertex \(i\), and let
\[
\widehat C_{\deg}=\Bigl\{i:\ \deg(i)>(1+\eta)V_m n r_n^m\bigl[(1-p_n)p_{00}^m+p_np_{10}^m\bigr]\Bigr\}
\]
and
\[
\widehat C_+=\widehat C_{\mathrm{star}}\cup\widehat C_{\deg}.
\]
The threshold is \((1+\eta)\) times the largest expected degree of a peripheral vertex, which is attained away from the boundary. Together with the stars, the degree rule then covers every position of the cube whenever
\begin{equation}\label{eq:Ceta}
2^{-m}(1-p_n)p_{10}^m\ \ge\ (1+\eta)^3\bigl[(1-p_n)p_{00}^m+p_np_{10}^m\bigr].
\end{equation}
This holds for large \(n\) and small \(\eta\) when
\[
p_n\to0,\qquad p_{10}>2p_{00}.
\]
In that case even a core vertex at a corner, whose ball is cut to a fraction \(2^{-m}\), has a larger mean degree than a peripheral vertex in the bulk.

Let
\[
c_1=\frac{V_m(p_{10}-p_{00})^m}{4^{m+1}},\qquad c_\eta=\frac{\eta^2V_mp_{00}^m}{16}\quad(\eta\in(0,1)),
\]
where \(c_1\) is the rate at which a core vertex fails to be the centre of an induced star, and \(c_\eta\) is the rate of the degree rule.

\begin{theorem}[Recovery of the core]\label{thm:recovery}
Under \(H_1\), the following hold.
\begin{enumerate}[label=\textup{(\roman*)},itemsep=2pt]
\item \textup{(Almost exact recovery.)} If \(np_n\to\infty\), \(r_n\to0\) and \(nr_n^m\to\infty\), then, for every \(\varepsilon>0\),
\[
\P\bigl(|\widehat C_{\mathrm{star}}\triangle C|\ge\varepsilon|C|\bigr)\to0 .
\]
\item \textup{(Exact recovery.)}
\begin{enumerate}[label=\textup{(\alph*)},itemsep=2pt]
\item If \(np_ne^{-c_1nr_n^m}\to0\) and \(np_nr_n\to0\), then
\[
\P(\widehat C_{\mathrm{star}}=C)\to1 .
\]
\item If \eqref{eq:Ceta} holds, \(np_ne^{-c_1nr_n^m}\to0\) and \(ne^{-c_\eta nr_n^m}\to0\), then
\[
\P(\widehat C_+=C)\to1 .
\]
\end{enumerate}
\end{enumerate}
\end{theorem}

Part~(i) requires three conditions. The expected degree must diverge, \(nr_n^m\to\infty\), so that a core vertex in the bulk is the centre of an induced star with high probability; the expected core size must diverge, \(np_n\to\infty\), so that the number of core vertices is large and a vanishing fraction of errors is meaningful; and \(r_n\to0\), so that the boundary layer of the cube, which has volume of order \(r_n\), becomes negligible. Under these conditions the star centres recover all but a vanishing fraction of the core, and the core vertices near the boundary, which are the hardest to find, are too few to matter.

Exact recovery requires more. Part~(ii)(a) shows that the star centres alone recover the whole core when two conditions hold. The first, \(np_ne^{-c_1nr_n^m}\to0\), holds whenever \(nr_n^m\) exceeds \(c_1^{-1}\log(np_n)\) by a diverging amount, for instance when \(nr_n^m\ge(1+\varepsilon)c_1^{-1}\log(np_n)\); it guarantees that every core vertex away from the boundary is the centre of an induced star. The second, \(np_nr_n\to0\), says that the expected number of core vertices near the boundary tends to zero; it is needed because the stars may miss such vertices. For \(m=1\), the two conditions are compatible: if \(nr_n=(1+\varepsilon)c_1^{-1}\log(np_n)\), then \(np_nr_n\le Cp_n\log n\), which tends to \(0\) as soon as \(p_n\log n\to0\). For \(m\ge2\), however, the boundary layer has volume of order \(r_n\gg r_n^m\), so the condition \(np_nr_n\to0\) is restrictive; this is where the degree rule is needed.

Part~(ii)(b) removes the boundary condition at the price of the stronger condition \(nr_n^m\gtrsim c_\eta^{-1}\log n\), which stems from the false positives of the degree rule: each of the \(n\) peripheral vertices has a small chance of having a large degree. Part~(ii), unlike part~(i), does not require \(np_n\to\infty\). When \(np_n\) is bounded, the core consists of a bounded number of vertices, which \(\widehat C_{\mathrm{star}}\) finds with high probability once the expected degree is large enough. Detection fails in that regime (Theorem~\ref{thm:main}(ii)) only because, with probability bounded away from zero, there is no core at all, and the graph is then homogeneous.

Since exact recovery implies almost exact recovery, Part~(ii)(a) also yields almost exact recovery when \(np_n=O(1)\), \(nr_n^m\to\infty\) and \(r_n\to0\). Combined with Part~(i), almost exact recovery is therefore possible whenever \(nr_n^m\to\infty\) and \(r_n\to0\), with no further restriction on \(np_n\).

The condition \(p_{10}>2p_{00}\) under which \eqref{eq:Ceta} holds is not an artefact of the proof: no rule based only on the degree can recover core vertices whose ball is cut too much by the boundary. A core vertex at such a position has a smaller expected degree than a typical peripheral vertex in the bulk, so any degree threshold that produces no false positives also misses it. When \(p_{10}<2p_{00}\), this includes the corners of the cube. The following proposition makes this precise.

\begin{proposition}[Limits of degree rules]\label{prop:degree}
Assume \(p_n\to0\) and \(nr_n^m\to\infty\), and let \(\eta\in(0,\frac18)\). Let \(\tau_n\) be non-random thresholds such that the rule \(\{i:\deg(i)>\tau_n\}\) contains no peripheral vertex with probability tending to one. Then, for every \(x\in[0,1]^m\) with \(\varphi(x)(p_{10}/p_{00})^m\le1-5\eta\), where \(\varphi(x)\in[2^{-m},1]\) is the fraction of the ball \(B(x,r_np_{10})\) that lies in the cube, a core vertex at \(x\) belongs to the rule with probability tending to \(0\).
\end{proposition}

\subsection{Sharpness of the recovery conditions}\label{sec:sharpness}

Let
\[
c_2=V_m(p_{11}^m-p_{00}^m)
\]
be the volume of the unit-ball shell between the radii \(p_{00}\) and \(p_{11}\); it governs the converse. An \emph{estimator} is any measurable function of the graph \(G\).

\begin{theorem}[Converse]\label{thm:sharp}
Assume~\eqref{eq:B}.
\begin{enumerate}[label=\textup{(\roman*)},itemsep=2pt]
\item \textup{(Exact recovery fails.)} Suppose that, for a fixed \(\varepsilon\in(0,1)\) and all large \(n\),
\[
nr_n^m\le(1-\varepsilon)\frac1{c_2}\log(np_n).
\]
Then, for every estimator \(\widehat C\),
\[
\liminf_{n\to\infty}\P\bigl(\widehat C\ne C\bigr)>0 .
\]
\item \textup{(Exact and almost exact recovery fail.)} Suppose \(nr_n^m\le K<\infty\) for all \(n\). Then, for every estimator \(\widehat C\),
\[
\liminf_{n\to\infty}\P\bigl(\widehat C\ne C\bigr)>0 ,
\qquad
\liminf_{n\to\infty}\P\bigl(|\widehat C\triangle C|>\delta|C|\bigr)>0,
\]
where \(\delta>0\) is a constant that depends only on \(K\), \(m\), \(p_{11}\), \(p_{10}\) and \(p_{00}\).
\end{enumerate}
\end{theorem}

Theorem~\ref{thm:sharp} says that, when the average degree is bounded, that is, \(nr_n^m\le K\), no estimator achieves exact or almost exact recovery. When the average degree grows more slowly than the order \(\log(np_n)\), no estimator achieves exact recovery.

Combining Theorems~\ref{thm:recovery} and~\ref{thm:sharp} yields sharp thresholds for almost exact recovery, together with conditions under which exact recovery is possible or impossible.

\begin{corollary}[Sharp and order-sharp recovery thresholds]\label{cor:sharp}
Assume \(r_n\to0\). Then the following hold.
\begin{enumerate}[label=\textup{(\alph*)},itemsep=2pt]
\item Almost exact recovery is possible if and only if \(nr_n^m\to\infty\).
\item Assume \(np_n\to\infty\). For any fixed \(\varepsilon>0\), the star centres \(\widehat C_{\mathrm{star}}\) achieve exact recovery if \(np_nr_n\to0\) and
\[
nr_n^m\ge(1+\varepsilon)\frac1{c_1}\log(np_n).
\]
\item Assume \(np_n\to\infty\). For any fixed \(\varepsilon>0\), the star centres together with the degree rule, \(\widehat C_+\), achieve exact recovery if \eqref{eq:Ceta} holds and
\[
nr_n^m\ge(1+\varepsilon)\max\Bigl\{\frac1{c_1}\log(np_n),\ \frac1{c_\eta}\log n\Bigr\}.
\]
\item Assume \(np_n\to\infty\). For any fixed \(\varepsilon>0\), exact recovery is impossible if
\[
nr_n^m\le(1-\varepsilon)\frac1{c_2}\log(np_n).
\]
\end{enumerate}
\end{corollary}

Corollary~\ref{cor:sharp} establishes that the threshold for almost exact recovery is sharp, and further implies that the threshold for exact recovery is sharp up to constant factors when the core has polynomial size. To see this, let \(p_n=n^{-t}\) for a fixed \(t\in(0,1)\), so that the expected number of core vertices \(np_n=n^{1-t}\) diverges and \(\log(np_n)=(1-t)\log n\) scales on the order of \(\log n\). Because the core size in real-world networks is typically much smaller than the total network size \(n\), the polynomial regime \(p_n=n^{-t}\) with a constant \(t\in(0,1)\) is a natural assumption that covers a wide range of CP structures \cite{YSM25}. By part~(d), exact recovery is impossible whenever
\[
nr_n^m\le(1-\varepsilon)c_2^{-1}(1-t)\log n .
\]
Conversely, by part~(c), combining the star centres with the degree rule achieves exact recovery when
\[
nr_n^m\ge(1+\varepsilon)\max\bigl\{(1-t)/c_1,1/c_\eta\bigr\}\log n,
\]
provided that~\eqref{eq:Ceta} holds, which requires \(p_{10}>2p_{00}\). Since both bounds scale as \(\log n\), the order of the threshold is sharp. Moreover, from the inequality
\[
(p_{10}-p_{00})^m\le(p_{11}-p_{00})^m\le p_{11}^m-p_{00}^m,
\]
it follows that \(c_2/c_1\ge4^{m+1}\), and thus
\[
(1-t)/c_2<(1-t)/c_1\le\max\bigl\{(1-t)/c_1,1/c_\eta\bigr\}.
\]
Consequently, the impossibility and achievability regimes are non-overlapping, confirming that the exact recovery threshold is sharp up to constant factors. Finally, as shown in part~(b), the star centres alone achieve the same \(\log n\) scaling without requiring \(p_{10}>2p_{00}\), provided that \(np_nr_n\to0\).

The recovery thresholds also clarify how recovery and detection compare. In block models, exact recovery implies almost exact recovery, which in turn implies detection, so the conditions for recovery are the stronger ones \cite{Abb18,MNS15}. The same ordering holds here when \(np_n\to\infty\). Detection needs only \(\liminf nr_n^m>0\) (Theorem~\ref{thm:main}(i)), whereas almost exact recovery needs \(nr_n^m\to\infty\) (Corollary~\ref{cor:sharp}(a)); moreover, recovery implies detection through the star test, since the event \(\widehat C_{\mathrm{star}}\ne\emptyset\) has no type~I error. In particular, if \(nr_n^m=O(1)\) and \(np_n\to\infty\), the core can be detected but not recovered almost exactly. The ordering is reversed only when \(np_n\) stays bounded: almost exact recovery then holds under \(nr_n^m\to\infty\) and \(r_n\to0\), while detection is impossible (Theorem~\ref{thm:main}(ii)). The reason is that the loss \(|\widehat C\triangle C|/|C|\) is measured relative to the size of the core, which may be empty or bounded. If \(C=\emptyset\), then \(\widehat C_{\mathrm{star}}=\emptyset\) is exactly right, because \(\widehat C_{\mathrm{star}}\subseteq C\) always; if \(C\ne\emptyset\), almost exact recovery of a bounded set is exact recovery, which the stars achieve once \(nr_n^m\to\infty\) and \(r_n\to0\). Hence, for bounded \(np_n\), the recovery statements describe the core \emph{given that it exists}; they do not say whether a core exists, which is what detection asks. In block models the communities have linear size, \(|C|\to\infty\), and this degenerate case does not occur.

For \(m=1\), the two constants can be made to meet. Theorem~\ref{thm:one} below improves the constant \(c_1\) of Corollary~\ref{cor:sharp}(b) to \(2(p_{10}-p_{00})\), which equals the converse constant \(c_2=2(p_{11}-p_{00})\) when \(p_{11}=p_{10}\). In this case, the exact recovery limit is sharp, including the constant, in the regime \(np_nr_n\to0\), and not only sharp in order.

\begin{theorem}[Sharp exact recovery for \(m=1\) and \(p_{11}=p_{10}\)]\label{thm:one}
Let \(m=1\), \(p_n\to0\) and \(np_n\to\infty\), and fix \(\varepsilon\in(0,1)\).
\begin{enumerate}[label=\textup{(\roman*)},itemsep=2pt]
\item The star centres achieve exact recovery if \(np_nr_n\to0\) and
\[
nr_n\ge(1+\varepsilon)\frac{\log(np_n)}{2(p_{10}-p_{00})} .
\]
\item No estimator achieves exact recovery if
\[
nr_n\le(1-\varepsilon)\frac{\log(np_n)}{2(p_{11}-p_{00})} .
\]
\end{enumerate}
In particular, if \(p_{11}=p_{10}\) and \(np_nr_n\to0\), the threshold for exact recovery of the core is sharp.
\end{theorem}

When \(p_{11}=p_{10}\), the CP-RGG is a random model for a noisy version of the ``ideal'' CP structure, in which the core vertices are connected to all other vertices and no two peripheral vertices are connected \cite{BE00,YSM25}. Therefore, for \(m=1\), Theorem~\ref{thm:one} gives sharp exact recovery limits for recovering a noisy ideal CP structure in the regime \(np_nr_n\to0\).

For \(m=1\), the constant \(1/c_2\) in the converse of Theorem~\ref{thm:sharp}(i) equals \(1/(2(p_{11}-p_{00}))\), and it agrees with the constant \(1/(2(p_{10}-p_{00}))\) in Theorem~\ref{thm:one}(i) when \(p_{11}=p_{10}\). The reason is geometric. A core vertex \(x\) is a star centre if three vertices in \([x-r_np_{10},x+r_np_{10}]\) are at mutual distances at least \(r_np_{00}\). On the line, the cheapest way for this to fail is that all peripheral vertices of the window lie in two intervals of length \(r_np_{00}\), which leaves a set of length \(2r_n(p_{10}-p_{00})\) empty; when \(p_{11}=p_{10}\), this is exactly the annulus around \(x\) that must be empty to hide its label. For \(m\ge2\), the constants in Corollary~\ref{cor:sharp}(b)--(d) do not match: the possibility constant comes from a convenient ball construction, and the impossibility constant from the empty annulus. The line argument does not extend, because \(s_m\) sets of diameter \(r_np_{00}\) have volume up to \(s_mV_m(r_np_{00}/2)^m>V_m(r_np_{00})^m\) (by the factor \(5/4\) for \(m=2\)). When \(p_{10}/p_{00}\) is large, the star test can therefore fail more often than the annulus is empty; its exact failure exponent is a packing problem that we leave open.

\section{Proofs}\label{sec:proofs}

\paragraph{Proof strategy.} All four theorems rest on one geometric fact and one symmetry. The geometric fact is a packing property: in a homogeneous random geometric graph a vertex never has \(s_m+1\) pairwise non-adjacent neighbours, whereas a core vertex has them with high probability, because its range \(r_np_{10}\) towards peripheral vertices exceeds their range \(r_np_{00}\) towards each other and leaves room for \(s_m+1\) well separated peripheral neighbours (Lemmas~\ref{lem:local}, \ref{lem:packing} and~\ref{lem:sparse}). The symmetry is that a core vertex with no other vertex in the annulus between the radii \(r_np_{00}\) and \(r_np_{11}\) has the same neighbours whatever its label.

For detection (Theorem~\ref{thm:main}), part~(i) uses the first fact: the test that looks for an induced star \(K_{1,s_m+1}\) never rejects the homogeneous null, and detects the alternative as soon as there are enough core vertices. Part~(ii) is a reduction. With probability \((1-p_n)^n\) the alternative has no core vertex and is then homogeneous (Lemma~\ref{lem:red}); the matched radius \(\gamma_n\) exceeds the periphery radius by a relative amount of order \(p_n\) (Lemma~\ref{lem:gamma}); and homogeneous graphs whose radii differ by a relative amount \(\delta\) cannot be perfectly separated when \(n\delta\) is bounded (Lemma~\ref{lem:scaling}).

For recovery (Theorem~\ref{thm:recovery}), the set \(\widehat C_{\mathrm{star}}\) of star centres is contained in \(C\) surely (Lemma~\ref{lem:centres}), so only missed core vertices count as errors. Each core vertex away from the boundary is missed with probability at most of order \(\exp(-c_1nr_n^m)\), and a union bound over the core, together with a concentration bound for the number of misses, gives almost exact and exact recovery; the degree rule then recovers the core vertices near the boundary. For the converse (Theorems~\ref{thm:sharp} and~\ref{thm:one}), the symmetry above hides the label of every core vertex whose annulus is empty. There are many such vertices, and no estimator can do better than guess their labels, which gives either a failure probability bounded below or a constant fraction of misclassified core vertices. Corollary~\ref{cor:sharp} follows by combining Theorems~\ref{thm:recovery} and~\ref{thm:sharp}, and Proposition~\ref{prop:degree} shows by a degree computation why the degree rule cannot recover core vertices at positions where the ball is cut too much by the boundary.

We use the same labels and positions \((Z,X)\) for all graphs. The constant \(c_1\) is the one defined before Theorem~\ref{thm:recovery}, namely
\[
c_1=\frac{V_m(p_{10}-p_{00})^m}{4^{m+1}} .
\]

\subsection{Lemmas}

\begin{lemma}[Edge probability]\label{lem:edgeprob}
Let \(X\) and \(Y\) be independent and uniform on \([0,1]^m\), let \(W=X-Y\), and for \(w=(w_1,\dots,w_m)\in\mathbb R^m\) let
\[
g(w)=\prod_{l=1}^m\bigl(1-|w_l|\bigr)_+ ,\qquad (t)_+=\max(t,0).
\]
Then \(g\) is a probability density of \(W\),
\[
f(\rho)=\P(\|W\|<\rho)=\int_{\|w\|<\rho}g(w)dw,
\]
and:
\begin{enumerate}[label=\textup{(\alph*)},itemsep=1pt]
\item \(f\) is continuous and strictly increasing on \([0,1]\), and \(f(\rho)\le V_m\rho^m\);
\item if \(\rho\le\rho'\le\frac14\), then
\[
f(\rho')-f(\rho)\ge\Bigl(\frac34\Bigr)^mV_m(\rho'^m-\rho^m);
\]
\item if \(\lambda\ge1\) and \(\lambda\rho\le1\), then \(f(\lambda\rho)\le\lambda^mf(\rho)\).
\end{enumerate}
\end{lemma}
\begin{proof}
\emph{Notation.} Here \(X=(X_1,\dots,X_m)\) and \(Y=(Y_1,\dots,Y_m)\) are two independent random points of the cube, the vector \(W=X-Y\in[-1,1]^m\) is the displacement between them, and \(w\in\mathbb R^m\) is the variable of integration, that is, a possible value of \(W\). By definition,
\[
f(\rho)=\P(\|X-Y\|<\rho)=\P(\|W\|<\rho).
\]

\emph{Step 1: \(g\) is the density of \(W\).} The \(2m\) coordinates \(X_1,\dots,X_m,Y_1,\dots,Y_m\) are independent and uniform on \([0,1]\). Hence the coordinates \(W_l=X_l-Y_l\), \(l=1,\dots,m\), are independent, because \(W_l\) is a function of \((X_l,Y_l)\) only, and each has the same density. For a single coordinate, the density of \(W_1=X_1-Y_1\) at \(t\in\mathbb R\) is \(h(t)=(1-|t|)_+\), the triangular density on \([-1,1]\). Since the coordinates of \(W\) are independent, the density of \(W\) is the product \(g(w)=\prod_{l=1}^mh(w_l)\). Consequently
\[
f(\rho)=\P(\|W\|<\rho)=\int_{\{\|w\|<\rho\}}g(w)dw .
\]
Two properties of \(g\) are used below: \(0\le g\le1\), and \(g(w)>0\) when \(\max_l|w_l|<1\); in particular \(g>0\) on the ball \(\{\|w\|<1\}\), because \(|w_l|\le\|w\|\).

\emph{(a) Monotonicity and the upper bound.} Since \(g\) is bounded, the map \(\rho\mapsto\int_{\{\|w\|<\rho\}}g\) is continuous, as the volume of the ball is continuous in \(\rho\). For \(0\le\rho<\rho'\le1\),
\[
f(\rho')-f(\rho)=\int_{\{\rho\le\|w\|<\rho'\}}g(w)dw>0,
\]
because the shell \(\{\rho\le\|w\|<\rho'\}\) has positive volume and \(g>0\) on it. This shows that \(f\) is strictly increasing. Using \(g\le1\),
\[
f(\rho)\le\operatorname{vol}\{\|w\|<\rho\}=V_m\rho^m.
\]

\emph{(b) A lower bound on increments.} If \(\|w\|<\rho'\le\frac14\), then \(|w_l|\le\|w\|<\frac14\) for every \(l\), so each factor of \(g\) is larger than \(\frac34\) and \(g(w)\ge(3/4)^m\). Integrating over the shell \(\{\rho\le\|w\|<\rho'\}\), whose volume is \(V_m(\rho'^m-\rho^m)\), gives
\[
f(\rho')-f(\rho)\ge\Bigl(\frac34\Bigr)^mV_m(\rho'^m-\rho^m).
\]

\emph{(c) Scaling.} Let \(\lambda\ge1\). Substitute \(w=\lambda v\), so that \(dw=\lambda^mdv\) and \(\|w\|<\lambda\rho\) if and only if \(\|v\|<\rho\):
\[
f(\lambda\rho)=\lambda^m\int_{\{\|v\|<\rho\}}g(\lambda v)dv\ \le\ \lambda^m\int_{\{\|v\|<\rho\}}g(v)dv=\lambda^mf(\rho),
\]
where the inequality holds because \((1-\lambda|v_l|)_+\le(1-|v_l|)_+\) for each \(l\) when \(\lambda\ge1\), so \(g(\lambda v)\le g(v)\).
\end{proof}

\begin{lemma}[Reduction to the case without a core]\label{lem:red}
Let \(\Omega=\{0,1\}^{\binom n2}\) be the set of graphs on the vertex set \(\{1,\dots,n\}\), let \(P_1\) denote the law on \(\Omega\) of the graph \(\mathcal G(r_n,p_n,p_{11},p_{10},p_{00})\), and let \(\mathcal G(r_np_{00})\) also denote the law of the homogeneous graph with radius \(r_np_{00}\). Then for every probability law \(Q\) on \(\Omega\),
\[
1-\TV(P_1,Q)\ \ge\ (1-p_n)^n\Bigl(1-\TV\bigl(\mathcal G(r_np_{00}),Q\bigr)\Bigr).
\]
\end{lemma}

The lemma says that the alternative cannot be told apart from \(Q\) better than the homogeneous graph with radius \(r_np_{00}\) can, up to the factor \((1-p_n)^n\), which is the probability that the alternative has no core vertex. The law \(Q\) is arbitrary: it need not be homogeneous. The lemma is used with \(Q=\mathcal G(\gamma_n)\) in the proof of Theorem~\ref{thm:main}(ii).

\begin{proof}
Write \(P_0\) for the law of \(\mathcal G(r_np_{00})\), \(c=(1-p_n)^n\in(0,1)\), and let \(E=\{Z_1=\dots=Z_n=0\}\) be the event that no vertex is a core vertex.

\emph{Step 1: \(P_1\) is a mixture.} Since \(Z_1,\dots,Z_n\) are i.i.d.\ Bernoulli\((p_n)\),
\[
\P(E)=(1-p_n)^n=c.
\]
On \(E\), every pair of vertices is a pair of peripheral vertices, so
\[
A_{ij}=\ind\{\|X_i-X_j\|<r_np_{00}\}
\]
for all \(i\ne j\). The labels are independent of the positions, so given \(E\) the positions are still i.i.d.\ uniform on the cube, and therefore the conditional law of the graph given \(E\) is \(P_0\). Let \(R\) be the conditional law of the graph given the complement \(E^c\). By the law of total probability, for every graph \(g\in\Omega\),
\[
P_1(g)=\P(A=g\mid E)\P(E)+\P(A=g\mid E^c)\P(E^c)=cP_0(g)+(1-c)R(g).
\]

\emph{Step 2: a bound on the overlap.} For \(g\in\Omega\), the minimum \(\min(\cdot,Q(g))\) is non-decreasing in its first argument, and \((1-c)R(g)\ge0\), so
\[
\min\bigl(P_1(g),Q(g)\bigr)\ \ge\ \min\bigl(cP_0(g),Q(g)\bigr).
\]
Moreover, since \(c\le1\) we have \(cQ(g)\le Q(g)\), hence
\[
\min\bigl(cP_0(g),Q(g)\bigr)\ \ge\ \min\bigl(cP_0(g),cQ(g)\bigr)=c\min\bigl(P_0(g),Q(g)\bigr).
\]

\emph{Step 3: sum over graphs.} For two laws \(P\) and \(Q\) on the finite set \(\Omega\), Definition~\ref{def:tv} gives
\[
1-\TV(P,Q)=\sum_{g\in\Omega}\min(P(g),Q(g)).
\]
Summing the two inequalities of Step 2 over \(g\in\Omega\),
\[
1-\TV(P_1,Q)=\sum_{g}\min\bigl(P_1(g),Q(g)\bigr)\ \ge\ c\sum_{g}\min\bigl(P_0(g),Q(g)\bigr)=c\bigl(1-\TV(P_0,Q)\bigr),
\]
which is the claim, since \(c=(1-p_n)^n\).
\end{proof}

\begin{lemma}[Scaling]\label{lem:scaling}
Let \(n\ge2\) and \(m\ge1\) be integers, and let \(\rho>0\) and \(\delta>0\). Let \(\mathcal G(\rho)\) and \(\mathcal G(\rho(1+\delta))\) denote the laws of the homogeneous random geometric graphs on \(n\) vertices with radii \(\rho\) and \(\rho(1+\delta)\) (the same formula as in Definition~\ref{def:hom} defines the graph for every radius \(\rho>0\)). Then
\[
\TV\bigl(\mathcal G(\rho),\mathcal G(\rho(1+\delta))\bigr)\ \le\ 1-(1+\delta)^{-mn}\ \le\ 1-e^{-mn\delta}.
\]
\end{lemma}

The lemma says that two homogeneous graphs whose radii differ by a relative amount \(\delta\) cannot be told apart perfectly, as long as \(mn\delta\) is bounded: a larger radius is the same as a smaller radius seen on a smaller cube, and a smaller cube is a sub-event of probability \((1+\delta)^{-mn}\). It is applied in Theorem~\ref{thm:main}(ii) with \(\rho=r_np_{00}\) and a \(\delta\) of order \(p_n\), for which \(mn\delta\) is of order \(np_n\).

\begin{proof}
\emph{Notation.} Let \(K=([0,1]^m)^n\) be the set of configurations \(y=(y_1,\dots,y_n)\) of \(n\) points of the cube, let \(U\) be the uniform law on \(K\), and for \(\rho'>0\) let \(\Phi_{\rho'}\) map a configuration \(y\) to the graph with adjacency matrix \(A_{ij}=\ind\{\|y_i-y_j\|<\rho'\}\), \(i\ne j\). If \(Y\sim U\), then \(\Phi_{\rho'}(Y)\) is the homogeneous graph with radius \(\rho'\). Put
\[
c=(1+\delta)^{-mn},\qquad S=\bigl([0,(1+\delta)^{-1}]^m\bigr)^n\subseteq K,
\]
so that \(c=\operatorname{vol}(S)=U(S)\), and write \(P=\mathcal G(\rho)\) and \(Q=\mathcal G(\rho(1+\delta))\).

\emph{Step 1: \(P\ge cQ\).} Let \(Y\sim U\). The event \(\{Y\in S\}\) has probability \(c\), and given this event \(Y\) is uniform on \(S\). Hence, with \(\Phi_\rho\) applied to \(Y\) and \(U_S\) the uniform law on \(S\),
\[
P(B)=\P\bigl(\Phi_\rho(Y)\in B\bigr)\ \ge\ \P\bigl(Y\in S\bigr)\P\bigl(\Phi_\rho(Y)\in B\mid Y\in S\bigr)=c\P\bigl(\Phi_\rho(Y_S)\in B\bigr),
\]
where \(Y_S\sim U_S\) and \(B\) is a set of graphs. It remains to identify the law of \(\Phi_\rho(Y_S)\).

\emph{Step 2: \(\Phi_\rho(Y_S)\sim Q\).} If \(Y'\sim U\), then \(Y'/(1+\delta)\), obtained by dividing every coordinate by \(1+\delta\), has independent coordinates that are uniform on \([0,(1+\delta)^{-1}]\), so \(Y'/(1+\delta)\sim U_S\). For every pair \(i\ne j\),
\[
\Bigl\|\frac{Y'_i}{1+\delta}-\frac{Y'_j}{1+\delta}\Bigr\|<\rho\iff\|Y'_i-Y'_j\|<\rho(1+\delta),
\]
so
\[
\Phi_\rho\Bigl(\frac{Y'}{1+\delta}\Bigr)=\Phi_{\rho(1+\delta)}(Y'),
\]
and the right-hand side has law \(Q\). Thus \(\Phi_\rho(Y_S)\sim Q\), and Step 1 reads \(P(B)\ge cQ(B)\) for every \(B\); since the set of graphs is finite, \(P(g)\ge cQ(g)\) for every graph \(g\).

\emph{Step 3: conclusion.} By Definition~\ref{def:tv},
\[
1-\TV(P,Q)=\sum_g\min(P(g),Q(g)).
\]
By Step 2,
\[
\min(P(g),Q(g))\ge\min(cQ(g),Q(g))=cQ(g),
\]
and summing over \(g\) gives
\[
1-\TV(P,Q)\ge c\sum_gQ(g)=c.
\]
Hence
\[
\TV(P,Q)\le1-(1+\delta)^{-mn}.
\]
Finally, \(1+\delta\le e^\delta\) gives \((1+\delta)^{-mn}\ge e^{-mn\delta}\).
\end{proof}

\begin{lemma}[Matched radius]\label{lem:gamma}
Assume~\eqref{eq:B}. Let \(\gamma_n\) be the density-matched radius of Section~\ref{sec:hyp}, and write \(\gamma_n=r_np_{00}(1+u)\), where \(u\ge0\) is the relative excess of \(\gamma_n\) over \(r_np_{00}\). Then
\[
u\ \le\ \frac{2}{m}\Bigl(\frac43\Bigr)^m(a-1)p_n,\qquad a=\Bigl(\frac{p_{11}}{p_{00}}\Bigr)^m .
\]
In particular \(mnu\le2(4/3)^m(a-1)np_n\).
\end{lemma}

The lemma says that the matched radius \(\gamma_n\) exceeds the peripheral radius \(r_np_{00}\) by a relative amount of order \(p_n\). This is natural: the mixture of the three radii differs from \(r_np_{00}\) only for the pairs that contain a core vertex, which form a fraction at most \(2p_n\) of all pairs. The bound is used in the proof of Theorem~\ref{thm:main}(ii), where it controls the factor \(e^{-mnu}\) from Lemma~\ref{lem:scaling}.

\begin{proof}
\emph{Step 0: setting.} Write \(f_{ab}=f(r_np_{ab})\) for \(p_{ab}\in\{p_{11},p_{10},p_{00}\}\), so that \(f_{00}=f(r_np_{00})\), and recall that
\[
f(\gamma_n)=p_n^2f_{11}+2p_n(1-p_n)f_{10}+(1-p_n)^2f_{00}
\]
by definition. The right-hand side is a weighted average of \(f_{11}\ge f_{10}\ge f_{00}\), since \(f\) is increasing, so it lies in \([f_{00},f_{11}]\); as \(f\) is continuous and strictly increasing (Lemma~\ref{lem:edgeprob}(a)), \(\gamma_n\) exists, is unique, and lies in \([r_np_{00},r_np_{11}]\). In particular \(u\ge0\), and \(\gamma_n\le r_np_{11}\le\frac14\) by~\eqref{eq:B}. Also \(f_{00}>0\), because \(f\) is strictly increasing with \(f(0)=0\). Put
\[
K=\frac{f(\gamma_n)}{f_{00}}-1,\qquad \pi=p_n^2+2p_n(1-p_n)=1-(1-p_n)^2\le2p_n .
\]
The proof bounds \(K\) from above in Step 1 and from below in Step 2, and compares the two bounds in Step 3.

\emph{Step 1: an upper bound on \(K\).} The three weights \(p_n^2\), \(2p_n(1-p_n)\) and \((1-p_n)^2\) add up to \(1\), so subtracting \(f_{00}\) from the definition of \(f(\gamma_n)\) leaves only the pairs that contain a core vertex:
\[
f(\gamma_n)-f_{00}=p_n^2\bigl(f_{11}-f_{00}\bigr)+2p_n(1-p_n)\bigl(f_{10}-f_{00}\bigr).
\]
Dividing by \(f_{00}\),
\[
K=p_n^2\Bigl(\frac{f_{11}}{f_{00}}-1\Bigr)+2p_n(1-p_n)\Bigl(\frac{f_{10}}{f_{00}}-1\Bigr).
\]
For \(p_{ab}\in\{p_{11},p_{10}\}\) apply Lemma~\ref{lem:edgeprob}(c) with \(\rho=r_np_{00}\) and \(\lambda=p_{ab}/p_{00}\ge1\), so that \(\lambda\rho=r_np_{ab}\le1\): it gives
\[
f(r_np_{ab})\le\Bigl(\frac{p_{ab}}{p_{00}}\Bigr)^mf_{00}\le af_{00},
\]
since \(p_{ab}\le p_{11}\). Hence each bracket is at most \(a-1\), and
\[
K\le\bigl(p_n^2+2p_n(1-p_n)\bigr)(a-1)=\pi(a-1)\le2p_n(a-1).
\]

\emph{Step 2: a lower bound on \(K\).} Since \(r_np_{00}\le\gamma_n\le\frac14\), Lemma~\ref{lem:edgeprob}(b) with \(\rho=r_np_{00}\) and \(\rho'=\gamma_n\) gives
\[
f(\gamma_n)-f_{00}\ge\Bigl(\frac34\Bigr)^mV_m\bigl(\gamma_n^m-(r_np_{00})^m\bigr),
\]
and Lemma~\ref{lem:edgeprob}(a) gives \(f_{00}\le V_m(r_np_{00})^m\). Therefore
\[
K=\frac{f(\gamma_n)-f_{00}}{f_{00}}\ \ge\ \frac{\Bigl(\frac34\Bigr)^mV_m\bigl(\gamma_n^m-(r_np_{00})^m\bigr)}{V_m(r_np_{00})^m}=\Bigl(\frac34\Bigr)^m\bigl((1+u)^m-1\bigr).
\]
By Bernoulli's inequality, \((1+u)^m\ge1+mu\), so
\[
K\ge(3/4)^mmu.
\]

\emph{Step 3: conclusion.} Combining Steps 1 and 2,
\[
(3/4)^mmu\le K\le2p_n(a-1),
\]
that is,
\[
u\le\frac2m(4/3)^m(a-1)p_n .
\]
Multiplying by \(mn\) gives
\[
mnu\le2(4/3)^m(a-1)np_n .
\]
\end{proof}

\begin{lemma}[Local packing]\label{lem:local}
Let \(x\in\mathbb R^m\), let \(\rho_1,\dots,\rho_k>0\), and let \(y_1,\dots,y_k\in\mathbb R^m\) be points such that
\[
\|y_j-x\|<\rho_j\quad\text{for all }j,\qquad \|y_j-y_{j'}\|\ge\max\{\rho_j,\rho_{j'}\}\quad\text{for all }j\ne j'.
\]
Then \(k\le s_m\), where \(s_m\) is the constant of Section~\ref{sec:notation}, the largest number of unit vectors of \(\mathbb R^m\) whose pairwise angles all exceed \(\pi/3\).
\end{lemma}

In the applications, \(x\) is the position of a vertex, \(y_1,\dots,y_k\) are the positions of \(k\) of its neighbours, and \(\rho_j\) is the radius at which the pair formed by the vertex and the \(j\)th neighbour is linked, so that \(\|y_j-x\|<\rho_j\). The second condition says that the neighbours are pairwise non-adjacent whenever the radius of the pair \(\{j,j'\}\) is at least \(\max\{\rho_j,\rho_{j'}\}\): the lemma bounds the number of pairwise non-adjacent neighbours of a vertex by \(s_m\), provided the radius of a pair of neighbours is not smaller than the radii of the two links to the vertex. This holds in a homogeneous graph, where all radii are equal, and for a peripheral vertex in the core--periphery graph (Lemmas~\ref{lem:packing} and~\ref{lem:centres}). The statement does not involve the cube or any random variable; it is a statement about points in \(\mathbb R^m\).

\begin{proof}
If \(k\le1\) there is nothing to prove, since \(s_m\ge1\); so let \(k\ge2\). The proof has three steps: the points \(y_j\) differ from \(x\), the directions from \(x\) to the \(y_j\) are well separated, and the number of such directions is at most \(s_m\).

\emph{Step 1: \(y_j\ne x\) for every \(j\).} Suppose \(y_j=x\) for some \(j\), and take \(j'\ne j\). Then
\[
\|y_j-y_{j'}\|=\|x-y_{j'}\|<\rho_{j'}\le\max\{\rho_j,\rho_{j'}\},
\]
which contradicts the assumption
\[
\|y_j-y_{j'}\|\ge\max\{\rho_j,\rho_{j'}\}.
\]
Hence the unit vectors \(u_j=(y_j-x)/\|y_j-x\|\) are well defined.

\emph{Step 2: the angles between the \(u_j\) exceed \(\pi/3\).} Fix \(j\ne j'\). Let \(M\) and \(\mu\) be the larger and the smaller of the two distances \(\|y_j-x\|\) and \(\|y_{j'}-x\|\), and let \(\theta\in[0,\pi]\) be the angle at \(x\) between \(y_j\) and \(y_{j'}\), that is, the angle between \(u_j\) and \(u_{j'}\). Suppose \(\theta\le\pi/3\). Then \(\cos\theta\ge\frac12\), and the law of cosines in the triangle with vertices \(x\), \(y_j\), \(y_{j'}\) gives
\[
\|y_j-y_{j'}\|^2=M^2+\mu^2-2M\mu\cos\theta\ \le\ M^2+\mu^2-M\mu=M^2-\mu(M-\mu)\ \le\ M^2 .
\]
Since \(\|y_j-x\|<\rho_j\) and \(\|y_{j'}-x\|<\rho_{j'}\), we have \(M<\max\{\rho_j,\rho_{j'}\}\). Hence \(\|y_j-y_{j'}\|<\max\{\rho_j,\rho_{j'}\}\), which contradicts the assumption. Therefore \(\theta>\pi/3\).

\emph{Step 3: counting.} By Step 2, \(u_1,\dots,u_k\) are \(k\) unit vectors of \(\mathbb R^m\) (distinct, since their pairwise angles are positive) whose pairwise angles all exceed \(\pi/3\). By the definition of \(s_m\), \(k\le s_m\).
\end{proof}

\begin{lemma}[Packing]\label{lem:packing}
Recall that an \emph{induced star} \(K_{1,k}\) in a graph consists of a vertex, the \emph{centre}, and \(k\) of its neighbours, the \emph{leaves}, no two of which are adjacent. Assume~\eqref{eq:B}.
\begin{enumerate}[label=\textup{(\alph*)},itemsep=1pt]
\item For every \(\rho>0\), the homogeneous graph \(\mathcal G(\rho)\) contains no induced star \(K_{1,s_m+1}\).
\item Under \(H_1\), the probability that the graph contains an induced star \(K_{1,s_m+1}\) is at least
\[
1-e^{-2^{-m}np_n}-e^{-n/8}-(s_m+1)\exp\Bigl(-c_1nr_n^m\Bigr).
\]
\end{enumerate}
\end{lemma}

Part~(a) says that a homogeneous graph never has a vertex with \(s_m+1\) pairwise non-adjacent neighbours, whatever the radius; this is why the star test has no type~I error. Part~(b) says that the core--periphery graph does contain such a vertex with high probability, namely a core vertex with \(s_m+1\) peripheral neighbours in well separated directions; the gap between the radii \(r_np_{10}\) and \(r_np_{00}\) is what leaves room for them.

\begin{proof}
\textbf{Part (a).} Let \(x\) be the position of the centre of an induced star \(K_{1,k}\) in \(\mathcal G(\rho)\), and \(y_1,\dots,y_k\) the positions of its leaves. Each leaf is adjacent to the centre, so \(\|y_j-x\|<\rho\). The leaves are pairwise non-adjacent, so \(\|y_j-y_{j'}\|\ge\rho\) for \(j\ne j'\). Lemma~\ref{lem:local} with all \(\rho_j=\rho\) gives \(k\le s_m\). Hence there is no induced \(K_{1,s_m+1}\).

\medskip
\textbf{Part (b).} The proof has four steps: a geometric construction of \(s_m+1\) separated balls around a point (Step~1), the observation that a core vertex with a peripheral vertex in each ball is the centre of an induced star (Step~2), the probability that such a core vertex exists (Step~3), and the conclusion (Step~4).

\emph{Step 1: the balls.} Put \(\rho_0=r_np_{00}\), \(\nu=p_{10}/p_{00}>1\), and
\[
R_0=\rho_0\frac{1+\nu}2,\qquad w=\frac{r_n(p_{10}-p_{00})}4=\rho_0\frac{\nu-1}4 .
\]
Let \(e_1,\dots,e_{s_m}\) be unit vectors with pairwise angles larger than \(\pi/3\), which exist by the definition of \(s_m\), and set \(\tau_0=0\) and \(\tau_k=R_0e_k\) for \(1\le k\le s_m\). Then \(\|\tau_0-\tau_k\|=R_0\), and for \(1\le k\ne k'\le s_m\) the angle \(\theta\) between \(e_k\) and \(e_{k'}\) exceeds \(\pi/3\), so
\[
\|\tau_k-\tau_{k'}\|=2R_0\sin\Bigl(\frac\theta2\Bigr)>R_0.
\]
Thus \(\|\tau_k-\tau_{k'}\|\ge R_0\) for all \(k\ne k'\). Let \(I=[r_np_{10},1-r_np_{10}]^m\) be the interior of the cube. Since \(r_np_{10}\le\frac14\) by~\eqref{eq:B},
\[
\P(X\in I)=(1-2r_np_{10})^m\ge2^{-m}.
\]
For \(x\in I\) and \(0\le k\le s_m\), let \(B_k(x)=B(x+\tau_k,w)\) be the open ball of radius \(w\) around \(x+\tau_k\). These \(s_m+1\) balls have the following properties.
\begin{enumerate}[label=\textup{(\roman*)},itemsep=1pt]
\item \emph{They lie in the ball \(B(x,r_np_{10})\), which lies in the cube.} For \(y\in B_k(x)\),
\[
\|y-x\|<R_0+w=\frac{\rho_0(3\nu+1)}{4}<\rho_0\nu=r_np_{10},
\]
because \(\nu>1\); and \(B(x,r_np_{10})\subseteq[0,1]^m\) because \(x\in I\). Hence each \(B_k(x)\) has volume \(V_mw^m\) inside the cube.
\item \emph{They are pairwise at distance more than \(\rho_0\).} For \(y\in B_k(x)\), \(y'\in B_{k'}(x)\) and \(k\ne k'\),
\[
\|y-y'\|>\|\tau_k-\tau_{k'}\|-2w\ge R_0-2w=\rho_0=r_np_{00}.
\]
In particular the balls are disjoint.
\end{enumerate}

\emph{Step 2: a core vertex with a peripheral vertex in each ball is a star centre.} Let \(i\) be a core vertex with \(X_i=x\in I\), and suppose that each ball \(B_k(x)\), \(0\le k\le s_m\), contains the position of a peripheral vertex, say \(j_k\). Then \(i\) is adjacent to every \(j_k\): the pair \((i,j_k)\) is a core--peripheral pair, with radius \(r_np_{10}\), and \(\|X_{j_k}-x\|<r_np_{10}\) by~(i). The vertices \(j_k\) are pairwise non-adjacent: each pair \((j_k,j_{k'})\) is peripheral, with radius \(r_np_{00}\), and \(\|X_{j_k}-X_{j_{k'}}\|>r_np_{00}\) by~(ii). Hence \(i\) is the centre of an induced \(K_{1,s_m+1}\) with leaves \(j_0,\dots,j_{s_m}\).

\emph{Step 3: such a core vertex exists with high probability.} Let \(M\) be the number of core vertices whose position lies in \(I\). Each vertex is a core vertex with position in \(I\) with probability \(p_n\P(X\in I)\ge2^{-m}p_n\), independently, so
\[
\P(M=0)=\bigl(1-p_n\P(X\in I)\bigr)^n\le e^{-2^{-m}np_n}.
\]
Let \(N\) be the number of core vertices. Then \(N\) is Binomial\((n,p_n)\) with mean at most \(n/2\) by~\eqref{eq:B}, and Hoeffding's inequality gives
\[
\P\Bigl(N>\frac{3n}{4}\Bigr)\le\P\Bigl(N-np_n\ge\frac n4\Bigr)\le\exp\Bigl(-2n\Bigl(\frac14\Bigr)^2\Bigr)=e^{-n/8}.
\]
It remains to bound the probability that, on the event \(\{M\ge1\}\cap\{N\le3n/4\}\), some ball around a suitable core vertex contains no peripheral vertex. This is done in four points.

\begin{enumerate}[label=\textup{(\arabic*)},itemsep=3pt]
\item \emph{Choice of the core vertex.} On the event \(\{M\ge1\}\) let \(i^*\) be the core vertex of smallest index whose position lies in \(I\), and write \(x=X_{i^*}\in I\). The vertex \(i^*\) is determined by the labels \(Z\) and by the positions of the core vertices only.
\item \emph{The peripheral positions are still uniform.} Condition on the labels \(Z\), restricted to the event \(\{N\le3n/4\}\), and on the positions of all core vertices, restricted to \(\{M\ge1\}\). The labels are independent of the positions, and the positions of different vertices are independent of each other. Hence, given \(Z\) and the positions of the core vertices, the positions of the peripheral vertices are still i.i.d.\ uniform on the cube, and they are independent of \(x\), which is a function of the conditioning variables. Their number is \(n-N\ge n/4\) on the event \(\{N\le3n/4\}\), and the vertex \(i^*\) is not among them since it is a core vertex.
\item \emph{One ball is empty with small probability.} Fix \(k\in\{0,\dots,s_m\}\). By property~(i) of Step~1, the ball \(B_k(x)\) lies in the cube, so a uniform point falls in it with probability exactly its volume \(V_mw^m\), and it avoids it with probability \(1-V_mw^m\). By independence of the \(n-N\) peripheral positions, and since \(1-t\le e^{-t}\) and \(n-N\ge n/4\), the probability that none of them lies in \(B_k(x)\) is
\[
\bigl(1-V_mw^m\bigr)^{n-N}\ \le\ \bigl(1-V_mw^m\bigr)^{n/4}\ \le\ \exp\Bigl(-\frac n4V_mw^m\Bigr).
\]
The exponent is, with \(w=r_n(p_{10}-p_{00})/4\),
\[
\frac n4V_mw^m=\frac n4V_m\Bigl(\frac{r_n(p_{10}-p_{00})}4\Bigr)^m=\frac{V_m(p_{10}-p_{00})^m}{4^{m+1}}nr_n^m=c_1nr_n^m,
\]
so the bound is \(\exp(-c_1nr_n^m)\).
\item \emph{Union bound over the balls.} There are \(s_m+1\) balls \(B_0(x),\dots,B_{s_m}(x)\), so the conditional probability that at least one of them contains no peripheral vertex is at most \((s_m+1)\exp(-c_1nr_n^m)\). This bound does not depend on the values of the conditioning variables, so it also holds after averaging over them.
\end{enumerate}

\emph{Step 4: conclusion.} If \(M\ge1\), \(N\le3n/4\) and every ball \(B_k(X_{i^*})\) contains a peripheral vertex, then \(i^*\) is the centre of an induced \(K_{1,s_m+1}\) by Step~2. The probability that one of these three conditions fails is at most
\[
e^{-2^{-m}np_n}+e^{-n/8}+(s_m+1)\exp\bigl(-c_1nr_n^m\bigr),
\]
by Step~3 and a union bound. This proves~(b).
\end{proof}

Lemma~\ref{lem:packing}(b) is informative only when the expected degree \(nr_n^m\) diverges, because its failure probability contains the term \(\exp(-c_1nr_n^m)\), which does not vanish when \(nr_n^m\) stays bounded. The next lemma covers this regime, which is needed for Theorem~\ref{thm:main}(i). When \(nr_n^m\le K\), a given core vertex in the bulk is the centre of an induced star only with small probability, of order \((nr_n^m)^{s_m+1}\), since each of the \(s_m+1\) regions in which a leaf may lie must contain a peripheral vertex. The argument therefore uses the large number of core vertices instead of the degree: among the roughly \(np_n\) core vertices one selects many that are well separated, so that their neighbourhoods are disjoint and the corresponding events are negatively associated, and the star exists with probability tending to one as soon as \(np_n(nr_n^m)^{s_m+1}\to\infty\). Under the standing assumption~\eqref{eq:B}, \(nr_n^m\) is bounded below, so this condition is implied by \(np_n\to\infty\).

\begin{lemma}[Stars at bounded expected degree]\label{lem:sparse}
Assume \(r_np_{11}\le\frac14\), \(p_n\le\frac12\), \(nr_n^m\le K\) for a fixed \(K<\infty\), and \(np_n(nr_n^m)^{s_m+1}\to\infty\). Then under \(H_1\) the graph contains an induced star \(K_{1,s_m+1}\) with probability tending to \(1\).
\end{lemma}

\begin{proof}
Write \(\lambda_n=nr_n^m\), \(\rho_0=r_np_{00}\), \(\rho_1=r_np_{10}\) and \(I=[\rho_1,1-\rho_1]^m\). Since \(\lambda_n\le K\), \(r_n\to0\), and since \(\lambda_n\le K\) and \(np_n\lambda_n^{s_m+1}\to\infty\), also \(np_n\to\infty\).

\emph{Step 1: the regions.} Choose unit vectors \(u_1,\dots,u_{s_m}\) with pairwise angles at least \(\alpha_*>\pi/3\) (possible by the definition of \(s_m\)), then \(\theta>0\) with \(2\theta<\alpha_*-\pi/3\), and \(\delta\in(0,\frac14)\) with \((1-2\delta)p_{10}>p_{00}\). For \(x\in I\) let
\[
\begin{aligned}
R_0(x)&=B(x,\delta\rho_1),\\
R_k(x)&=\bigl\{y:(1-\delta)\rho_1<\|y-x\|<\rho_1,\ \angle(y-x,u_k)<\theta\bigr\},\quad 1\le k\le s_m .
\end{aligned}
\]
These are disjoint subsets of \(B(x,\rho_1)\subseteq[0,1]^m\), and each has volume at least \(v_*r_n^m\) for a constant \(v_*>0\) depending only on \(m\), \(p_{10}\), \(p_{00}\) (the volume of \(R_k\) is \(\omega V_m\rho_1^m(1-(1-\delta)^m)\), where \(\omega\in(0,1)\) is the fraction of the unit ball in a cone of half-angle \(\theta\)). Let \(y_k\in R_k(x)\), \(0\le k\le s_m\). Then \(\|y_k-x\|<\rho_1\) for all \(k\), and \(\|y_k-y_l\|>\rho_0\) for \(k\ne l\): for \(k,l\ge1\), the angle between \(y_k-x\) and \(y_l-x\) exceeds \(\pi/3\), so with \(a=\|y_k-x\|\) and \(b=\|y_l-x\|\),
\[
\|y_k-y_l\|^2>a^2+b^2-ab\ge\min(a,b)^2\ge(1-\delta)^2\rho_1^2>\rho_0^2 ;
\]
and for \(k\ge1\), \(l=0\),
\[
\|y_k-y_0\|\ge(1-\delta)\rho_1-\delta\rho_1=(1-2\delta)\rho_1>\rho_0.
\]
Consequently, if a core vertex \(i\) with \(X_i=x\in I\) has a peripheral vertex in each \(R_k(x)\), \(0\le k\le s_m\), then \(i\) and these \(s_m+1\) peripheral vertices form an induced star \(K_{1,s_m+1}\): the leaves are adjacent to \(i\) because the core--peripheral radius is \(\rho_1\), and pairwise non-adjacent because the peripheral radius is \(\rho_0\). Call this event \(E_i\).

\emph{Step 2: one core vertex.} Condition on the labels and on the positions of the core vertices, and suppose \(N=|C|\le3n/4\), which has probability at least \(1-e^{-n/8}\) by Chernoff's bound since \(p_n\le\frac12\). The \(n'=n-N\ge n/4\) peripheral positions are i.i.d. uniform. Split the peripheral vertices into \(s_m+1\) groups of size at least \(n/(8(s_m+1))\) for large \(n\). The events \(\{\text{group }k\text{ has a vertex in }R_k(x)\}\) are independent and imply \(E_i\), so for \(X_i=x\in I\),
\[
\P(E_i)\ \ge\ \Bigl(1-\bigl(1-v_*r_n^m\bigr)^{n/(8(s_m+1))}\Bigr)^{s_m+1}\ \ge\ \bigl(1-e^{-c_5\lambda_n}\bigr)^{s_m+1},
\]
and by \(1-e^{-t}\ge\min(t,1)/2\) and \(\lambda_n\le K\),
\[
\P(E_i)\ \ge\ \Bigl(\frac{c_6\lambda_n}{2}\Bigr)^{s_m+1}=:q_n,
\]
with \(c_5=v_*/(8(s_m+1))\) and \(c_6=\min(c_5,1/K)\). Note that \(np_nq_n\to\infty\) by assumption.

\emph{Step 3: many separated core vertices.} Condition on the labels, hence on \(N=|C|\); the positions \(X_i\), \(i\in C\), are then i.i.d. uniform on \([0,1]^m\). Let \(M\) be the number of core vertices with position in \(I\), and let \(W\) be the number of ordered pairs \((i,j)\) of distinct core vertices \(i,j\) with \(X_i,X_j\in I\) and \(\|X_i-X_j\|<2\rho_1\). We show that, with probability at least \(1-2v_n-1/L\) given \(N\in[\frac12np_n,2np_n]\), there are at least \(N/(2(1+2Lc_3))\) core vertices in \(I\) whose mutual distances are at least \(2\rho_1\).

\emph{(a) Most core vertices lie in \(I\).} The complement \([0,1]^m\setminus I\) is contained in the union of \(2m\) slabs of width \(\rho_1\) (one at each face), so
\[
\P(X\notin I)\le2m\rho_1=:v_n\to0 .
\]
Hence \(\E[N-M\mid N]\le v_nN\), and by Markov's inequality
\[
\P\Bigl(M<\frac N2\mid N\Bigr)=\P\Bigl(N-M>\frac N2\mid N\Bigr)\le2v_n .
\]

\emph{(b) Few close pairs.} For distinct core vertices \(i,j\), conditionally on \(X_i\) the position \(X_j\) is uniform, so
\[
\P(\|X_i-X_j\|<2\rho_1)\le V_m(2\rho_1)^m.
\]
There are fewer than \(N^2\) ordered pairs, and \(\rho_1^m=p_{10}^mr_n^m\), so
\[
\E[W\mid N]\le N^2V_m(2\rho_1)^m=N^22^mV_mp_{10}^mr_n^m .
\]
If \(N\le2np_n\), then \(Nr_n^m\le2np_nr_n^m=2p_n\lambda_n\), hence
\[
\E[W\mid N]\le N\cdot2^{m+1}V_mp_{10}^mp_n\lambda_n\le c_3N,\qquad c_3=2^mV_mp_{10}^mK,
\]
because \(p_n\le\frac12\) and \(\lambda_n\le K\). For \(L\ge1\), Markov's inequality gives
\[
\P\bigl(W>Lc_3N\mid N\bigr)\le\frac1L .
\]

\emph{(c) An independent set.} Consider the graph \(\Gamma\) whose vertices are the \(M\) core vertices with position in \(I\), two of them being joined if their distance is less than \(2\rho_1\). Each edge of \(\Gamma\) is counted twice in \(W\), so \(\Gamma\) has \(W/2\) edges and average degree \(\bar d=W/M\). By the Caro--Wei bound and Jensen's inequality (the function \(d\mapsto1/(d+1)\) is convex), \(\Gamma\) has an independent set \(T\) with
\[
|T|\ \ge\ \sum_{v}\frac1{d_v+1}\ \ge\ \frac M{1+\bar d}\ =\ \frac{M^2}{M+W}.
\]
The function \(M\mapsto M^2/(M+W)\) is increasing, so on the event \(\{M\ge N/2,\ W\le Lc_3N\}\),
\[
|T|\ \ge\ \frac{(N/2)^2}{N/2+Lc_3N}\ =\ \frac N{2(1+2Lc_3)} .
\]

\emph{(d) Disjoint regions.} If \(i\ne j\) are in \(T\), then \(\|X_i-X_j\|\ge2\rho_1\), so the open balls \(B(X_i,\rho_1)\) and \(B(X_j,\rho_1)\) are disjoint. Since \(R_k(X_i)\subseteq B(X_i,\rho_1)\), all the regions \(R_k(X_i)\), \(i\in T\), \(0\le k\le s_m\), are pairwise disjoint (regions of the same vertex are disjoint by Step~1).

\emph{Step 4: negative association.} Let \(\mathcal F\) denote the core data, that is, the labels \(Z_1,\dots,Z_n\) and the positions \(X_i\), \(i\in C\). Work on the event \(\{N\le3n/4\}\). Given \(\mathcal F\), the \(n'=n-N\) peripheral positions are i.i.d.\ uniform on the cube, and the set \(T\) of Step~3(c) is a fixed set of core vertices, because it is a deterministic function of the core positions. We show that
\[
\P\Bigl(\bigcap_{i\in T}E_i^c\ \Big|\ \mathcal F\Bigr)\ \le\ \prod_{i\in T}\bigl(1-\P(E_i\mid\mathcal F)\bigr)\ \le\ (1-q_n)^{|T|}.
\]
The events \(E_i\), \(i\in T\), are not independent, since there are exactly \(n'\) peripheral vertices and a vertex lying in one region cannot lie in another. The dependence is, however, of the favourable sign, and this is what we prove.

\begin{enumerate}[label=\textup{(\alph*)},itemsep=2pt]
\item \emph{The counts are multinomial.} Let \(N_{i,k}\) be the number of peripheral vertices in \(R_k(X_i)\), \(i\in T\), \(0\le k\le s_m\), and let \(N_{\mathrm{rem}}\) be the number in the remainder of the cube. By Step~3(d) the regions \(R_k(X_i)\) are pairwise disjoint, so together with the remainder they partition the cube. Since the \(n'\) peripheral positions are i.i.d.\ uniform, the vector \((N_{i,k},N_{\mathrm{rem}})\) is, given \(\mathcal F\), multinomial with \(n'\) trials and cell probabilities equal to the volumes of the cells. The cells are determined by \(\mathcal F\), which is why we condition on it.
\item \emph{Negative association.} A random vector \(\xi\) is negatively associated if \(\E[f(\xi_A)g(\xi_B)]\le\E[f(\xi_A)]\E[g(\xi_B)]\) for all disjoint index sets \(A,B\) and all nonnegative functions \(f,g\) that are both nondecreasing or both nonincreasing in each coordinate. Joag-Dev and Proschan \cite{JDP83} showed that multinomial vectors are negatively associated.
\item \emph{The events are decreasing in disjoint families of counts.} We have
\[
\ind_{E_i^c}=\ind\bigl\{N_{i,k}=0\ \text{for some }k\bigr\}=h_i\bigl(N_{i,0},\dots,N_{i,s_m}\bigr),
\]
where \(h_i\in\{0,1\}\) is nonincreasing in each coordinate: increasing a count can only turn ``some region is empty'' into ``no region is empty''. The function \(h_i\) depends only on the counts of the regions of \(i\), and these families of counts are disjoint for different \(i\in T\).
\item \emph{Product bound.} Enumerate \(T=\{i_1,\dots,i_t\}\) and put \(H_j=\prod_{l\le j}h_{i_l}\). Each \(H_j\) is nonnegative and nonincreasing in each coordinate, and depends only on the counts of \(i_1,\dots,i_j\), while \(h_{i_{j+1}}\) depends on a disjoint set of counts. By~(b),
\[
\E\bigl[H_{j+1}\mid\mathcal F\bigr]\le\E\bigl[H_j\mid\mathcal F\bigr]\E\bigl[h_{i_{j+1}}\mid\mathcal F\bigr],
\]
and induction on \(j\) gives
\[
\P\Bigl(\bigcap_{i\in T}E_i^c\ \Big|\ \mathcal F\Bigr)=\E\bigl[H_t\mid\mathcal F\bigr]\le\prod_{i\in T}\P\bigl(E_i^c\mid\mathcal F\bigr).
\]
\item \emph{Conclusion.} Every \(i\in T\) has \(X_i\in I\), so by Step~2, \(\P(E_i\mid\mathcal F)\ge q_n\). Hence each factor is at most \(1-q_n\), and there are \(|T|\) factors. This proves the claim.
\end{enumerate}

\emph{Step 5: conclusion.} As \(np_n\to\infty\), Chebyshev's inequality gives
\[
\P\bigl(N\notin[\tfrac12np_n,2np_n]\bigr)\le\frac{4}{np_n}\to0.
\]
Combining Steps~2--4, for each fixed \(L\ge1\),
\[
\P_1(\text{no induced star})\le e^{-n/8}+\frac4{np_n}+2v_n+\frac1L
+\exp\Bigl(-\frac{q_nnp_n}{4(1+2Lc_3)}\Bigr),
\]
and the right-hand side has \(\limsup\le1/L\), because \(np_nq_n\to\infty\). Since \(L\) is arbitrary, the probability tends to \(0\).
\end{proof}

\subsection{Proof of Theorem~\ref{thm:main}}

\paragraph{(i)} The null \(\mathcal G(\gamma_n)\) is homogeneous, so by Lemma~\ref{lem:packing}(a) the star test has no type~I error. Let \(np_n\to\infty\). Since
\[
\TV\ge\P_1(\text{reject})-\P_0(\text{reject}),
\]
it suffices to show that every subsequence has a further subsequence along which the type~II error tends to \(0\). Along a subsequence, either \(nr_n^m\to\infty\) or \(nr_n^m\le K\) for some \(K<\infty\) along a further subsequence.

\emph{Case \(nr_n^m\to\infty\).} By Lemma~\ref{lem:packing}(b) the type~II error is at most
\[
e^{-2^{-m}np_n}+e^{-n/8}+(s_m+1)\exp\bigl(-c_1nr_n^m\bigr)\to0 .
\]

\emph{Case \(nr_n^m\le K\).} By~\eqref{eq:B}, \(nr_n^m\ge\lambda_0\) for some \(\lambda_0>0\) and large \(n\), so
\[
np_n(nr_n^m)^{s_m+1}\ge np_n\lambda_0^{s_m+1}\to\infty.
\]
By Lemma~\ref{lem:sparse} the type~II error tends to \(0\).

In both cases \(\TV\to1\) along the subsequence, hence along the whole sequence.

\paragraph{(ii)} By Lemma~\ref{lem:red} with \(Q=\mathcal G(\gamma_n)\) and Lemma~\ref{lem:scaling} with \(\rho=r_np_{00}\) and \(\delta=u\) from Lemma~\ref{lem:gamma},
\begin{equation}
1-\TV(H_0,H_1)\ge(1-p_n)^ne^{-mnu}\ge e^{-2np_n}\exp\Bigl(-2\Bigl(\frac43\Bigr)^m(a-1)np_n\Bigr)=e^{-C_*np_n}\label{eq:scal}
\end{equation}
for every \(n\), where \(a=(p_{11}/p_{00})^m\), using \((1-p_n)^n\ge e^{-2np_n}\) (since \(\log(1-p_n)\ge-p_n/(1-p_n)\ge-2p_n\) for \(p_n\le\frac12\)) and \(mnu\le2(4/3)^m(a-1)np_n\) by Lemma~\ref{lem:gamma}. If \(np_n\le K\), then \(1-\TV\ge e^{-C_*K}\), which is the claim.

\paragraph{The equivalence.} If \(np_n\to\infty\), part~(i) gives \(\TV\to1\). If \(np_n\not\to\infty\), then \(np_n\le C\) along a subsequence and part~(ii) gives \(\TV\le1-e^{-C_*C}\) along it, so \(\TV\not\to1\). \hfill\(\blacksquare\)

\subsection{Proof of Theorem~\ref{thm:recovery}}

The proof uses a geometric lemma on star centres, which gives the absence of false positives (F1), and a probability bound for a single vertex, from which the two parts of the theorem follow.

\begin{lemma}[Star centres are core vertices]\label{lem:centres}
Under \(H_1\), every vertex \(i\) that is the centre of an induced \(K_{1,s_m+1}\) is a core vertex.
\end{lemma}

The lemma says that, in the core--periphery graph, a peripheral vertex never has \(s_m+1\) pairwise non-adjacent neighbours. For a homogeneous graph this is Lemma~\ref{lem:packing}(a); here the radii of the pairs of neighbours differ, and the point is that they are never smaller than the radii of the links to the vertex.

\begin{proof}
Let \(i\) be a peripheral vertex and let \(j_1,\dots,j_k\) be pairwise non-adjacent neighbours of \(i\). We check the hypotheses of Lemma~\ref{lem:local}. Put \(x=X_i\) and \(y_h=X_{j_h}\), and let \(\rho_h=r_np_{10}\) if \(j_h\) is a core vertex and \(\rho_h=r_np_{00}\) if \(j_h\) is peripheral; these are the radii of the pairs \((i,j_h)\). Since \(i\) and \(j_h\) are adjacent, \(\|y_h-x\|<\rho_h\). For \(h\ne h'\), the radius of the pair \(\{j_h,j_{h'}\}\) is \(r_np_{00}\) if both are peripheral, \(r_np_{10}\) if exactly one is a core vertex, and \(r_np_{11}\) if both are core vertices.
In each case this radius is at least \(\max\{\rho_h,\rho_{h'}\}\): in the first case \(\rho_h=\rho_{h'}=r_np_{00}\); in the second case the larger of \(\rho_h,\rho_{h'}\) is \(r_np_{10}\); in the third case both equal \(r_np_{10}\le r_np_{11}\). Since \(j_h\) and \(j_{h'}\) are not adjacent, their distance is at least the radius of the pair, so
\[
\|y_h-y_{h'}\|\ge\max\{\rho_h,\rho_{h'}\}.
\]
Lemma~\ref{lem:local} gives \(k\le s_m\). Hence a peripheral vertex is not the centre of an induced \(K_{1,s_m+1}\).
\end{proof}

\paragraph{Notation and two preliminary bounds.} Let \(C=\{i:Z_i=1\}\) be the core, \(N=|C|\), \(I=[r_np_{10},1-r_np_{10}]^m\), and recall that \(\widehat C_{\mathrm{star}}\) is the set of centres of induced stars \(K_{1,s_m+1}\) and \(c_1\) and \(c_\eta\) are as in Theorem~\ref{thm:recovery}. For a vertex \(i\), let \(E_i\) be the event that \(i\) is a core vertex with \(X_i\in I\) that is \emph{missed}, that is, \(i\notin\widehat C_{\mathrm{star}}\). We use two bounds.

\begin{enumerate}[label=\textup{(P\arabic*)},itemsep=2pt,leftmargin=*]
\item \(\P\bigl(N>\tfrac{3n}{4}\bigr)\le e^{-n/8}\). This is Hoeffding's inequality, as in Step~3 of the proof of Lemma~\ref{lem:packing}(b), since \(N\) is Binomial\((n,p_n)\) with \(p_n\le\frac12\).
\item For every \(i\),
\[
\P\bigl(E_i\cap\{N\le\tfrac{3n}{4}\}\bigr)\le p_n(s_m+1)\exp\bigl(-c_1nr_n^m\bigr).
\]
\end{enumerate}

\emph{Proof of (P2).} Fix \(i\). Condition on the whole vector of labels \(Z\), restricted to the event \(\{Z_i=1,\ N\le3n/4\}\), and on the position \(x=X_i\), restricted to \(x\in I\). The positions are independent of the labels, so the positions of the \(n-N\ge n/4\) peripheral vertices are still i.i.d.\ uniform on the cube, independent of \(x\). By Steps~1 and~2 of the proof of Lemma~\ref{lem:packing}(b), if each of the \(s_m+1\) balls \(B_k(x)\) contains the position of a peripheral vertex, then \(i\) is the centre of an induced \(K_{1,s_m+1}\), so \(i\in\widehat C_{\mathrm{star}}\). Each ball has volume \(V_mw^m\) inside the cube, with \(w=r_n(p_{10}-p_{00})/4\), so it contains none of the peripheral positions with probability at most
\[
\bigl(1-V_mw^m\bigr)^{n/4}\le\exp\bigl(-\frac n4V_mw^m\bigr)=\exp(-c_1nr_n^m).
\]
A union bound over the \(s_m+1\) balls shows that, given \(Z\) and \(x\) as above, \(i\) is missed with probability at most \((s_m+1)\exp(-c_1nr_n^m)\). Averaging over \(Z\) and \(x\), and using \(\P(Z_i=1)=p_n\), gives (P2). \(\square\)

\paragraph{The degree rule.} Write \(\rho_0=r_np_{00}\), \(\rho_1=r_np_{10}\),
\[
\lambda_0=V_mr_n^m\bigl[(1-p_n)p_{00}^m+p_np_{10}^m\bigr],
\]
\(\tau=(1+\eta)n\lambda_0\), so that
\[
\widehat C_{\deg}=\{i:\deg(i)>\tau\}.
\]
Note \(\lambda_0\ge V_m\rho_0^m\). For \(x\in[0,1]^m\), let \(\varphi(x)\in[2^{-m},1]\) be the fraction of the ball \(B(x,\rho_1)\) that lies in the cube (it equals \(1\) on \(I\) and \(2^{-m}\) at a corner), and put
\[
D_\eta=\Bigl\{x\in[0,1]^m:\ (1-p_n)\varphi(x)p_{10}^m\ \ge\ (1+\eta)^3\bigl[(1-p_n)p_{00}^m+p_np_{10}^m\bigr]\Bigr\},
\]
the set of positions at which a core vertex has expected degree at least \((1+\eta)^3(n-1)\lambda_0\) (counting only its peripheral neighbours), which exceeds \((1+\eta)\tau\) for \(n\ge2/\eta\), so that its degree exceeds \(\tau\) with high probability, and
\[
u_n=\operatorname{vol}\bigl([0,1]^m\setminus(I\cup D_\eta)\bigr)\le2m\rho_1,
\]
the volume of the positions covered by neither rule. Since \(\varphi\ge2^{-m}\), we have \(u_n=0\) as soon as~\eqref{eq:Ceta} holds. The degree rule is analysed in two steps.

\emph{Step D1: the law of a degree.} Condition on \(X_i=x\) and \(Z_i\). The pairs \((Z_j,X_j)\), \(j\ne i\), are i.i.d., so the indicators \(\ind\{j\sim i\}\) are i.i.d.\ Bernoulli and \(\deg(i)\) is Binomial\((n-1,q)\) with
\[
q=p_n\operatorname{vol}\bigl(B(x,r_np_{Z_i1})\cap[0,1]^m\bigr)+(1-p_n)\operatorname{vol}\bigl(B(x,r_np_{Z_i0})\cap[0,1]^m\bigr),
\]
where \(B(x,\rho)\) is the Euclidean ball.

\emph{Step D2: the two bounds.} (a) \emph{A peripheral vertex.} Here \(Z_i=0\), and since a ball meets the cube in at most its volume,
\[
q\le V_mr_n^m\bigl[(1-p_n)p_{00}^m+p_np_{10}^m\bigr]=\lambda_0.
\]
By monotonicity of the binomial law in \(q\) and the Chernoff bound
\[
\P(\mathrm{Bin}\ge(1+\delta)\mu)\le e^{-\delta^2\mu/3}
\]
for \(\delta\in(0,1]\), with \(\mu=(n-1)\lambda_0\) and \(\tau\ge(1+\eta)\mu\),
\[
\P(\deg(i)>\tau)\le\exp\Bigl(-\frac{\eta^2(n-1)\lambda_0}{3}\Bigr)\le e^{-c_\eta nr_n^m},
\]
because \((n-1)\ge n/2\) and \(\lambda_0\ge V_mp_{00}^mr_n^m\).

(b) \emph{A core vertex with \(X_i=x\in D_\eta\).} Here \(Z_i=1\), and counting only the peripheral neighbours,
\[
q\ge(1-p_n)\operatorname{vol}(B(x,\rho_1)\cap[0,1]^m)=(1-p_n)\varphi(x)V_m\rho_1^m\ge(1+\eta)^3\lambda_0
\]
by the definition of \(D_\eta\). For \(n\ge2/\eta\), \((n-1)/n\ge(1+\eta)^{-1}\), so the mean \(\mu_1=(n-1)q\) satisfies
\[
\mu_1\ge(1+\eta)^2n\lambda_0=(1+\eta)\tau.
\]
Consequently
\[
\tau\le\frac{\mu_1}{1+\eta}=(1-\delta)\mu_1,\qquad\text{where}\quad\delta=1-\frac1{1+\eta}=\frac{\eta}{1+\eta}.
\]
Since \(\eta\in(0,1)\), we have \(1+\eta<2\), and hence
\[
\delta\in(0,1),\qquad \delta>\frac\eta2 .
\]
Because \(\tau\le(1-\delta)\mu_1\), the event \(\{\deg(i)\le\tau\}\) is contained in \(\{\deg(i)\le(1-\delta)\mu_1\}\). Given \(X_i=x\) and \(Z_i=1\), \(\deg(i)\) is binomial with \(n-1\) trials and success probability \(q\), so the lower-tail bound
\[
\P(\mathrm{Bin}\le(1-\delta)\mu_1)\le e^{-\delta^2\mu_1/2}
\]
gives
\[
\P(\deg(i)\le\tau)\le e^{-\delta^2\mu_1/2}\le e^{-\eta^2\mu_1/8}\le e^{-c_\eta nr_n^m},
\]
where the second inequality uses \(\delta>\eta/2\), and the last uses
\[
\mu_1\ge(n-1)\lambda_0\ge \frac{nV_mp_{00}^mr_n^m}{2}.
\]

\paragraph{No false positives, (F1).} If \(i\in\widehat C_{\mathrm{star}}\), then \(i\) is the centre of an induced \(K_{1,s_m+1}\), so \(i\) is a core vertex by Lemma~\ref{lem:centres}. Hence \(\widehat C_{\mathrm{star}}\subseteq C\) for every realisation, and in particular \(C\triangle\widehat C_{\mathrm{star}}=C\setminus\widehat C_{\mathrm{star}}\).

\paragraph{No false positives of the degree rule, (F2).} By Step~D2(a) and a union bound over the at most \(n\) peripheral vertices,
\[
\P(\widehat C_{\deg}\not\subseteq C)\le ne^{-c_\eta nr_n^m}\qquad\Bigl(n\ge\frac2\eta\Bigr),
\]
which tends to \(0\) under the condition \(ne^{-c_\eta nr_n^m}\to0\) of Theorem~\ref{thm:recovery}(ii)(b).

\paragraph{Proof of (i).} The proof has three steps: a bound on the expected number of missed core vertices, a lower bound on the size of the core, and the passage to convergence in probability.

\emph{Step 1: the expected number of misses.} A missed core vertex \(i\in C\setminus\widehat C_{\mathrm{star}}\) is of one of two kinds: its position lies outside \(I\), or it lies in \(I\) and \(E_i\) occurs. For the first kind, the labels and the positions are independent, so the expected number is \(np_n\P(X\notin I)\), and by Bernoulli's inequality \((1-t)^m\ge1-mt\),
\[
\P(X\notin I)=1-(1-2r_np_{10})^m\le2mr_np_{10}.
\]
For the second kind, on the event \(\{N>3n/4\}\) the number of misses is at most \(N\le n\), and on the event \(\{N\le3n/4\}\) we use (P2). Hence the expected number of misses of the second kind is at most
\[
n\P\bigl(N>\tfrac{3n}{4}\bigr)+\sum_{i=1}^n\P\bigl(E_i\cap\{N\le\tfrac{3n}{4}\}\bigr)\ \le\ ne^{-n/8}+np_n(s_m+1)\exp\bigl(-c_1nr_n^m\bigr).
\]
Adding the two bounds gives, for every \(n\),
\[
\E|C\setminus\widehat C_{\mathrm{star}}|\le ne^{-n/8}+np_n\Bigl[2mr_np_{10}+(s_m+1)\exp\bigl(-c_1nr_n^m\bigr)\Bigr].
\]

\emph{Step 2: the relative error in expectation.} Assume \(np_n\to\infty\), \(r_n\to0\) and \(nr_n^m\to\infty\). Dividing the inequality of Step~1 by \(np_n\),
\[
\frac{\E|C\setminus\widehat C_{\mathrm{star}}|}{np_n}\ \le\ \frac{e^{-n/8}}{p_n}+2mr_np_{10}+(s_m+1)\exp\bigl(-c_1nr_n^m\bigr).
\]
The last two terms tend to \(0\), because \(r_n\to0\) and \(nr_n^m\to\infty\). For the first term, \(np_n\ge1\) for large \(n\), so \(1/p_n\le n\) and \(e^{-n/8}/p_n\le ne^{-n/8}\to0\). Hence
\[
\frac{\E|C\setminus\widehat C_{\mathrm{star}}|}{np_n}\to0.
\]

\emph{Step 3: convergence in probability.} The number \(N\) of core vertices has mean \(np_n\) and variance \(np_n(1-p_n)\le np_n\), so by Chebyshev's inequality
\[
\P\Bigl(N<\frac{np_n}2\Bigr)\le\P\Bigl(|N-np_n|>\frac{np_n}2\Bigr)\le\frac{np_n}{\bigl(\frac{np_n}{2}\bigr)^2}=\frac{4}{np_n}\to0 .
\]
Let \(\varepsilon>0\). On the event \(\{N\ge np_n/2\}\), \(|C\setminus\widehat C_{\mathrm{star}}|>\varepsilon N\) implies \(|C\setminus\widehat C_{\mathrm{star}}|>\varepsilon np_n/2\). By Markov's inequality,
\[
\begin{aligned}
&\P\bigl(|C\setminus\widehat C_{\mathrm{star}}|>\varepsilon N\bigr)\ \le\ \P\Bigl(N<\frac{np_n}2\Bigr)+\P\Bigl(|C\setminus\widehat C_{\mathrm{star}}|>\frac{\varepsilon np_n}2\Bigr)\\
&\le\ \frac{4}{np_n}+\frac{2\E|C\setminus\widehat C_{\mathrm{star}}|}{\varepsilon np_n}\ \longrightarrow\ 0 .
\end{aligned}
\]
By (F1),
\[
\frac{|C\triangle\widehat C_{\mathrm{star}}|}{|C|}=\frac{|C\setminus\widehat C_{\mathrm{star}}|}{N},
\]
so \(|C\triangle\widehat C_{\mathrm{star}}|/|C|\to0\) in probability. \hfill\(\blacksquare\)

\paragraph{Proof of (ii)(a).} By (F1), \(\widehat C_{\mathrm{star}}\subseteq C\), so \(\{\widehat C_{\mathrm{star}}\ne C\}=\{C\not\subseteq\widehat C_{\mathrm{star}}\}\). A core vertex is missed either because \(N>3n/4\), or because its position lies outside \(I\), or because it lies in \(I\) and \(E_i\) occurs. The labels and positions are independent, so the expected number of core vertices outside \(I\) is \(np_n\P(X\notin I)\), and by Bernoulli's inequality \((1-t)^m\ge1-mt\),
\[
\P(X\notin I)=1-(1-2r_np_{10})^m\le2mr_np_{10}.
\]
Hence, by the union bound, (P1) and (P2),
\[
\begin{aligned}
&\P(\widehat C_{\mathrm{star}}\ne C)\le\P\bigl(N>\tfrac{3n}{4}\bigr)+\sum_{i=1}^n\P\bigl(E_i\cap\{N\le\tfrac{3n}{4}\}\bigr)+np_n\P(X\notin I)\\
&\le e^{-n/8}+np_n(s_m+1)\exp\bigl(-c_1nr_n^m\bigr)+2mr_np_{10}np_n ,
\end{aligned}
\]
for every \(n\). If \(np_ne^{-c_1nr_n^m}\to0\) and \(np_nr_n\to0\), the bound tends to \(0\), because \(e^{-n/8}\to0\) and \(s_m\) is fixed, so
\[
\P(\widehat C_{\mathrm{star}}=C)\to1.
\]

\paragraph{Proof of (ii)(b).} Let \(n\ge2/\eta\). Define the events
\begin{align*}
A_1&=\{\widehat C_{\deg}\not\subseteq C\},\\
A_2&=\{N>3n/4\}\cup\bigcup_{i=1}^n\bigl(E_i\cap\{N\le3n/4\}\bigr),\\
A_3&=\{\exists i\in C:\ X_i\in D_\eta\setminus I,\ i\notin\widehat C_{\deg}\},\\
A_4&=\{\exists i\in C:\ X_i\in[0,1]^m\setminus(I\cup D_\eta)\}.
\end{align*}
We first show that
\[
\{\widehat C_+\ne C\}\subseteq A_1\cup A_2\cup A_3\cup A_4 .
\]
Suppose that \(\widehat C_+\ne C\) and \(A_1\) does not occur. Then \(\widehat C_{\deg}\subseteq C\), and \(\widehat C_{\mathrm{star}}\subseteq C\) by~(F1), so \(\widehat C_+\subseteq C\). Hence some core vertex \(i\) belongs to neither \(\widehat C_{\mathrm{star}}\) nor \(\widehat C_{\deg}\). Its position lies in exactly one of \(I\), \(D_\eta\setminus I\) and \([0,1]^m\setminus(I\cup D_\eta)\).
\begin{itemize}[itemsep=2pt]
\item If \(X_i\in I\), then \(i\notin\widehat C_{\mathrm{star}}\) means that \(E_i\) occurs, so \(A_2\) occurs, whether or not \(N>3n/4\).
\item If \(X_i\in D_\eta\setminus I\), then \(i\notin\widehat C_{\deg}\) gives \(A_3\).
\item If \(X_i\notin I\cup D_\eta\), then \(A_4\) occurs.
\end{itemize}
We bound the four events.
\begin{itemize}[itemsep=2pt]
\item \emph{\(A_1\): a false positive of the degree rule.} By~(F2),
\[
\P(A_1)\le ne^{-c_\eta nr_n^m}.
\]
\item \emph{\(A_2\): \(N\) is large, or the stars miss a core vertex with position in \(I\).} By the union bound, (P1) and (P2), as in~(ii)(a),
\[
\P(A_2)\le e^{-n/8}+np_n(s_m+1)e^{-c_1nr_n^m}.
\]
\item \emph{\(A_3\): the degree rule misses a core vertex with position in \(D_\eta\setminus I\).} The label and the position of a vertex are independent, and the conditional law of \(\deg(i)\) given \((X_i,Z_i)\) is described in Step~D1. Hence, by Step~D2(b), for each \(i\),
\[
\P\bigl(Z_i=1,\ X_i\in D_\eta\setminus I,\ \deg(i)\le\tau\bigr)\le p_n\P(X\notin I)e^{-c_\eta nr_n^m},
\]
because \(\operatorname{vol}(D_\eta\setminus I)\le\P(X\notin I)\). Summing over \(i\), and using \(\P(X\notin I)\le1\),
\[
\P(A_3)\le np_n\P(X\notin I)e^{-c_\eta nr_n^m}\le np_ne^{-c_\eta nr_n^m}.
\]
\item \emph{\(A_4\): a core vertex in the region covered by neither rule.} A vertex is a core vertex with \(X_i\notin I\cup D_\eta\) with probability \(p_nu_n\), so
\[
\P(A_4)\le np_nu_n .
\]
\end{itemize}
Adding the four bounds gives, for \(n\ge2/\eta\),
\[
\P(\widehat C_+\ne C)\le e^{-n/8}+np_n(s_m+1)e^{-c_1nr_n^m}+n(1+p_n)e^{-c_\eta nr_n^m}+np_nu_n .
\]
If \eqref{eq:Ceta} holds, then \(u_n=0\), and if the two exponential conditions hold, every term tends to \(0\), because \(n(1+p_n)e^{-c_\eta nr_n^m}\le2ne^{-c_\eta nr_n^m}\). Finally, if \(nr_n^m\ge(1+\varepsilon)\log n/c_\eta\), then
\[
ne^{-c_\eta nr_n^m}\le n^{-\varepsilon}\to0.
\]
Note that the boundary term \(2mr_np_{10}np_n\) of part~(ii)(a) does not appear: the core vertices near the boundary are now found by the degree rule, which is why part~(ii)(b) does not require \(np_nr_n\to0\).

\subsection{Proof of Theorem~\ref{thm:sharp}}

\paragraph{Idea.} A core vertex \(i\) whose position \(X_i\) has no other vertex in the annulus between the radii \(r_np_{00}\) and \(r_np_{11}\) has the same neighbours whether it is a core vertex or a peripheral one, because every radius lies between \(r_np_{00}\) and \(r_np_{11}\). Its label is then independent of the graph, so no estimator can do better than guess it, and a guess is wrong with probability at least \(p_n\). The proof makes this precise (Step~1), shows that there are many such vertices (Step~2), and concludes (parts~(i) and~(ii)).

\paragraph{Notation.} Let \(\rho_0=r_np_{00}\) and \(\rho_2=r_np_{11}\), let \(J=[\frac14,\frac34]^m\), and for \(x\in\mathbb R^m\) let
\[
\mathrm{Ann}(x)=\{y:\rho_0\le\|y-x\|<\rho_2\}.
\]
Put
\[
c_2=V_m\bigl(p_{11}^m-p_{00}^m\bigr),\qquad a_1=V_m\bigl(\rho_2^m-\rho_0^m\bigr)=c_2r_n^m .
\]
The hypothesis of part~(i) reads \(nr_n^m\le(1-\varepsilon)\log(np_n)/c_2\). Both hypotheses of the theorem give \(nr_n^m=O(\log n)\): in part~(i),
\[
nr_n^m\le c_2^{-1}\log(np_n)\le c_2^{-1}\log n,
\]
and in part~(ii), \(nr_n^m\le K\). Hence \(r_n\to0\) and \(a_1\to0\), and for large \(n\) we have \(\rho_2<\frac14\). Then for \(x\in J\) the ball \(B(x,\rho_2)\) lies in \([0,1]^m\), so for a uniform point \(Y\) of the cube,
\[
\P(Y\in\mathrm{Ann}(x))=a_1.
\]

\paragraph{Step 1: hidden labels.} Let \(F\) be the set of vertices \(i\) with \(X_i\in J\) such that no other vertex lies in \(\mathrm{Ann}(X_i)\). We show that the labels in \(F\) cannot be learned from the graph, in three points.
\begin{enumerate}[label=\textup{(\arabic*)},itemsep=2pt]
\item \emph{The edges at \(F\) do not depend on the labels.} Let \(i\in F\) and let \(j\ne i\). Then either \(\|X_j-X_i\|<\rho_0\) or \(\|X_j-X_i\|\ge\rho_2\). In the first case \(i\) and \(j\) are adjacent, because every pair has radius at least \(\rho_0\); in the second case they are not adjacent, because every pair has radius at most \(\rho_2\) and adjacency requires distance strictly smaller than the radius. In both cases the answer does not depend on \(Z_i\) or \(Z_j\). Hence the graph \(G\) is a function of \((X,Z_{F^c})\), where \(Z_{F^c}\) are the labels outside \(F\).
\item \emph{The labels in \(F\) are still i.i.d.\ Bernoulli.} The set \(F\) is a function of \(X\) only, and the labels are i.i.d.\ and independent of \(X\). Therefore, conditionally on \(X\) and \(Z_{F^c}\), the labels \(Z_F\) are i.i.d.\ Bernoulli\((p_n)\), and they are independent of \(G\).
\item \emph{Consequences.} Condition on \((X,Z_{F^c})\). Then \(F\), \(G\), and hence \(\widehat C\), are fixed, while the labels \(Z_i\), \(i\in F\), are independent Bernoulli\((p_n)\) by point~(2). For \(i\in F\), the estimator's guess of \(Z_i\) is the constant \(g_i=\mathbf 1\{i\in\widehat C\}\in\{0,1\}\), which is determined by the graph and carries no information about \(Z_i\). If \(g_i=1\), the guess is correct with probability \(\P(Z_i=1)=p_n\) and wrong with probability \(1-p_n\); if \(g_i=0\), it is correct with probability \(1-p_n\) and wrong with probability \(p_n\). Since \(p_n\le\frac12\), we have \(p_n\le1-p_n\), so in both cases the guess is correct with probability at most \(1-p_n\) and wrong with probability at least \(p_n\). The events \(\{Z_i=g_i\}\) are independent over \(i\in F\). Hence, conditionally on \(X\) (so that \(|F|\) is fixed),
\[
\P\bigl(\text{all labels in }F\text{ are guessed correctly}\bigr)\le(1-p_n)^{|F|}\le e^{-p_n|F|},
\]
and the number of wrong guesses in \(F\) stochastically dominates a Binomial\((|F|,p_n)\) variable.
\end{enumerate}

\paragraph{Step 2: \(F\) is large.} Let \(S=|F|\). We compute its mean, then bound its variance and apply Chebyshev's inequality.

\emph{The mean.} A vertex has \(X_i\in J\) with probability \(2^{-m}\), and given \(X_i=x\in J\) each of the other \(n-1\) vertices avoids \(\mathrm{Ann}(x)\) independently with probability \(1-a_1\). Hence
\[
\E S=2^{-m}n(1-a_1)^{n-1}.
\]

\emph{The second moment.} We write
\[
\operatorname{Var}S=\E[S(S-1)]+\E S-(\E S)^2
\]
and bound
\[
\E[S(S-1)]=\sum_{i\ne j}\P(i,j\in F),
\]
splitting the ordered pairs according to the distance \(\|X_i-X_j\|\).
\begin{itemize}[itemsep=2pt]
\item \emph{Far pairs, \(\|X_i-X_j\|\ge2\rho_2\).} The annuli \(\mathrm{Ann}(X_i)\) and \(\mathrm{Ann}(X_j)\) lie in the disjoint balls \(B(X_i,\rho_2)\) and \(B(X_j,\rho_2)\), and neither vertex lies in the other's annulus. Hence the conditional probability that both belong to \(F\) is
\[
(1-2a_1)^{n-2}\le(1-a_1)^{2n-4},
\]
and these pairs contribute at most
\[
n^2 2^{-2m}(1-a_1)^{2n-4}=(1-a_1)^{-2}(\E S)^2.
\]
\item \emph{Near pairs, \(\|X_i-X_j\|<2\rho_2\).} Condition on \(i\in F\) and \(X_i=x\). The other \(n-1\) vertices are i.i.d.\ uniform on \([0,1]^m\setminus\mathrm{Ann}(x)\), whose density is at most \((1-a_1)^{-1}\le2\) for large \(n\), so the expected number of them within distance \(2\rho_2\) of \(x\) is at most \(2nV_m(2\rho_2)^m=C_0nr_n^m\), with \(C_0=2^{m+1}V_mp_{11}^m\). Hence these pairs contribute at most \(C_0nr_n^m\E S\) to \(\E[S(S-1)]\).
\end{itemize}
Therefore
\[
\operatorname{Var}S\ \le\ (1+C_0nr_n^m)\E S+\bigl((1-a_1)^{-2}-1\bigr)(\E S)^2 ,
\]
and Chebyshev's inequality gives
\[
\P\Bigl(S\le\frac{\E S}2\Bigr)\ \le\ \frac{4\operatorname{Var}S}{(\E S)^2}\ \le\ \frac{4(1+C_0nr_n^m)}{\E S}+4\bigl((1-a_1)^{-2}-1\bigr).
\]
The last term tends to \(0\) because \(a_1\to0\).

\paragraph{Proof of part (i).} Here we assume \(np_n\to\infty\); the general case is treated at the end of this part. The proof shows that \(\E S\) is large and then combines Step~1 with Step~2.

\emph{Step (i-a): the mean of \(S\).} Since \(\log(1-a)\ge-a-a^2\) for \(0\le a\le\frac12\), we have
\[
(1-a_1)^{n-1}\ge(1-a_1)^n\ge e^{-na_1(1+a_1)}.
\]
By hypothesis
\[
na_1=c_2nr_n^m\le(1-\varepsilon)\log(np_n).
\]
Since \(a_1\to0\), for large \(n\) we have
\[
(1-\varepsilon)(1+a_1)\le1-\frac\varepsilon2,
\]
and therefore
\[
\E S\ \ge\ 2^{-m}n(np_n)^{-(1-\varepsilon)(1+a_1)}\ \ge\ 2^{-m}n(np_n)^{-1+\varepsilon/2}=2^{-m}\frac{(np_n)^{\varepsilon/2}}{p_n}.
\]
Consequently
\[
p_n\E S\ge2^{-m}(np_n)^{\varepsilon/2}\to\infty,
\]
because \(np_n\to\infty\) and \(\varepsilon>0\), so that \((np_n)^{\varepsilon/2}\to\infty\).
Moreover, since \(nr_n^m\le c_2^{-1}\log(np_n)\) and \(p_n\le1\),
\[
\frac{\E S}{nr_n^m}\ \ge\ \frac{2^{-m}c_2(np_n)^{\varepsilon/2}}{p_n\log(np_n)}\ \ge\ \frac{2^{-m}c_2(np_n)^{\varepsilon/2}}{\log(np_n)}\ \longrightarrow\ \infty ,
\]
because \(np_n\to\infty\). Hence \((1+C_0nr_n^m)/\E S\to0\), and Step~2 gives
\[
\P(S\le\E S/2)\to0.
\]

\emph{Step (i-b): conclusion.} The event \(\{\widehat C=C\}\) requires all labels to be guessed correctly, in particular those in \(F\). By Step~1, its probability is at most \(\E[e^{-p_nS}]\), and
\[
\E\bigl[e^{-p_nS}\bigr]\ \le\ \P\Bigl(S\le\frac{\E S}2\Bigr)+e^{-p_n\E S/2}\ \longrightarrow\ 0 ,
\]
by Step (i-a). \hfill\(\square\)

\paragraph{Proof of part (i), general case.} The case \(np_n\to\infty\) was proved above. For general \(np_n\ge1\), let \(\ell=\liminf_n\P(\widehat C\ne C)\in[0,1]\), which exists because the probabilities lie in \([0,1]\). By the definition of the \(\liminf\), there is a subsequence along which \(\P(\widehat C\ne C)\to\ell\), and we pass to a further subsequence, along which the limit is still \(\ell\), such that either \(np_n\to\infty\) or \(np_n\le M\) for some \(M<\infty\). We show that \(\ell>0\) in both cases.

\emph{Case \(np_n\to\infty\).} The argument above uses only the hypothesis of part~(i), \(a_1\to0\) and \(np_n\to\infty\), so it applies along the subsequence and gives \(\P(\widehat C=C)\to0\). Hence \(\ell=1\).

\emph{Case \(np_n\le M\).} Since \(r_n>0\), the hypothesis of part~(i) forces \(\log(np_n)>0\), that is, \(np_n>1\), and it gives \(nr_n^m\le K:=c_2^{-1}\log M\). Hence part~(ii) gives \(\ell>0\).

In both cases \(\ell>0\), that is, \(\liminf\P(\widehat C\ne C)>0\). \hfill\(\blacksquare\)

\paragraph{Proof of part (ii), case \(np_n\to\infty\).} Here \(nr_n^m\le K\), so \(na_1=c_2nr_n^m\le c_2K\) and \(a_1\to0\). The proof has three steps.

\emph{Step (ii-a): \(S\) is of order \(n\).} As in part (i),
\[
(1-a_1)^{n-1}\ge e^{-na_1(1+a_1)}\ge e^{-2c_2K}
\]
for large \(n\), so \(\E S\ge c_Kn\) with \(c_K=2^{-m}e^{-2c_2K}\). Since \(nr_n^m\le K\), \((1+C_0nr_n^m)/\E S\to0\), and Step~2 gives
\[
\P(S\ge c_Kn/2)\to1.
\]

\emph{Step (ii-b): many errors in \(F\).} On the event \(\{S\ge c_Kn/2\}\), by Step~1 the number of wrong guesses in \(F\) dominates a Binomial\((\lceil c_Kn/2\rceil,p_n)\) variable \(B\), whose mean is at least \(c_Knp_n/2\to\infty\). By Chernoff's bound,
\[
\P(B\le\E B/2)\le e^{-\E B/8}\to0.
\]
Hence, with probability tending to one, the number of vertices of \(F\) with a wrong label is at least \(c_Knp_n/4\). Each such vertex belongs to \(\widehat C\triangle C\), so \(|\widehat C\triangle C|\ge c_Knp_n/4\).

\emph{Step (ii-c): comparison with \(|C|\).} The number of core vertices \(N=|C|\) has mean \(np_n\) and variance at most \(np_n\), so by Chebyshev's inequality \(\P(N>2np_n)\le1/(np_n)\to0\). Hence with probability tending to one,
\[
|\widehat C\triangle C|\ \ge\ \frac{c_Knp_n}4\ \ge\ \frac{c_K}8N,
\]
and \(\delta=c_K/8\), which depends only on \(K\), \(m\), \(p_{11}\) and \(p_{00}\), works. \hfill\(\square\)

\paragraph{Proof of part (ii), case \(np_n\) bounded.} Suppose \(1\le np_n\le M\). Steps~1 and~2 use only \(nr_n^m\le K\) and \(p_n\le\frac12\), not \(np_n\to\infty\); in particular
\[
\E S\ge c_Kn,\qquad \P(S<c_Kn/2)\to0
\]
as in Step (ii-a). By Step~1, conditionally on \(X\),
\[
\P(\widehat C=C\mid X)\le(1-p_n)^{S}\le e^{-p_nS},
\]
so
\[
\P(\widehat C=C)\le\P\bigl(S<\tfrac{c_Kn}2\bigr)+e^{-c_Knp_n/2}\le o(1)+e^{-c_K/2},
\]
which is the first claim of part~(ii) when \(np_n\) is bounded. For the second claim, on \(|\widehat C\triangle C|\), let \(E=\{\widehat C\ne C\}\), \(N=|C|\), and let \(\delta>0\) satisfy \(\delta M\le\frac12(1-e^{-c_K/2})\). On \(E\cap\{N<1/\delta\}\) we have \(|\widehat C\triangle C|\ge1>\delta N\) (also when \(N=0\)), so
\[
\P\bigl(|\widehat C\triangle C|>\delta N\bigr)\ge\P\bigl(E\cap\{N<1/\delta\}\bigr).
\]
By Markov's inequality,
\[
\P\Bigl(N\ge\frac1\delta\Bigr)\le\delta\E N=\delta np_n\le\delta M.
\]
Recall that \(\P(E_1\cap E_2)\ge\P(E_1)-\P(E_2^c)\) for any two events \(E_1\) and \(E_2\). Applying this with \(E_1=E\) and \(E_2=\{N<1/\delta\}\) in the previous display, we obtain
\[
\P\bigl(|\widehat C\triangle C|>\delta N\bigr)\ge\P(E)-\delta M\ge1-e^{-c_K/2}-o(1)-\delta M\ge\tfrac12\bigl(1-e^{-c_K/2}\bigr)-o(1),
\]
\hfill\(\square\)

\paragraph{Proof of part (ii), general case.} The first claim of part~(ii), \(\liminf\P(\widehat C\ne C)>0\), was shown in the bounded case with no upper bound on \(np_n\), because only \(np_n\ge1\) was used. The second claim, on \(|\widehat C\triangle C|\), was shown for \(np_n\to\infty\) in Steps~(ii-a)--(ii-c). It remains to prove it for general \(np_n\ge1\). For this, fix
\[
M_0=\max\Bigl(8,\ \frac{16\log8}{c_K}\Bigr),\qquad
\delta=\min\Bigl(\frac{c_K}{16},\ \frac{1-e^{-c_K/2}}{2M_0}\Bigr),\qquad
c=\frac12\bigl(1-e^{-c_K/2}\bigr),
\]
and note that \(c<\frac12\). We show that for every \(n\), whatever \(np_n\ge1\),
\[
\P\bigl(|\widehat C\triangle C|>\delta|C|\bigr)\ \ge\ c-\varepsilon_n,\qquad\varepsilon_n:=\P\bigl(S<\tfrac{c_Kn}{2}\bigr)\to0,
\]
where \(\varepsilon_n\) depends only on \(n\), \(K\), \(m\) and the radii ratios, not on \(p_n\) (Step~(ii-a)).
\begin{itemize}
\item If \(np_n>M_0\): on \(\{S\ge c_Kn/2\}\), the proof of Steps~(ii-b) and~(ii-c), carried out for fixed \(n\), shows that the event \(\{|\widehat C\triangle C|\ge\frac{c_K}8N\}\) fails with probability at most \(e^{-c_Knp_n/16}+1/(np_n)\) (Chernoff's bound for \(B\) with \(\E B\ge c_Knp_n/2\), and Chebyshev's bound for \(N\)). Since \(np_n>M_0\), this is at most
\[
e^{-c_KM_0/16}+1/M_0\le\frac18+\frac18.
\]
As \(\delta<c_K/8\) and
\[
|\widehat C\triangle C|\ge c_Knp_n/4>0,
\]
we get
\[
\P(|\widehat C\triangle C|>\delta N)\ge\frac34-\varepsilon_n\ge c-\varepsilon_n.
\]
\item If \(1\le np_n\le M_0\): the argument of the bounded case with \(M=M_0\) gives
\[
\P(|\widehat C\triangle C|>\delta N)\ge1-e^{-c_K/2}-\varepsilon_n-\delta M_0\ge c-\varepsilon_n,
\]
because \(\delta M_0\le\frac12(1-e^{-c_K/2})\).
\end{itemize}
Taking \(\liminf\) over \(n\) gives the claim. \hfill\(\blacksquare\)

\subsection{Proof of Theorem~\ref{thm:one}}

\paragraph{Part (i).} Let \(m=1\), \(\rho_0=r_np_{00}\), \(\rho_1=r_np_{10}\) and \(I=[\rho_1,1-\rho_1]\). The proof shows that a core vertex in \(I\) that is missed by \(\widehat C_{\mathrm{star}}\) must have all the peripheral points of its window in two short intervals (Steps~1 and~2), bounds the probability of this for one vertex (Step~3), and takes a union bound, adding the core vertices outside \(I\) (Step~4). By Lemma~\ref{lem:centres}, \(\widehat C_{\mathrm{star}}\subseteq C\), so only core vertices can be missed.

\emph{Step 1: a sufficient condition for a core vertex to be a star centre.} Let \(i\) be a core vertex with position \(x\in I\), and let \(W=(x-\rho_1,x+\rho_1)\subseteq[0,1]\) be its window. A peripheral vertex at a point of \(W\) is adjacent to \(i\), because a core--peripheral pair has radius \(\rho_1\). Suppose there are peripheral points \(a<b<c\) in \(W\) with \(b-a\ge\rho_0\) and \(c-b\ge\rho_0\). Then \(c-a\ge2\rho_0\) as well, so the three points are pairwise non-adjacent, since a peripheral pair has radius \(\rho_0\) and adjacency needs distance smaller than the radius. So \(i\) is the centre of an induced \(K_{1,3}\), and \(s_1+1=3\), hence \(i\in\widehat C_{\mathrm{star}}\). Thus if \(i\notin\widehat C_{\mathrm{star}}\), no such triple exists.

\emph{Step 2: if no triple exists, the peripheral points of \(W\) lie in two intervals of length \(\rho_0\).} Let \(a\) be the leftmost peripheral point of \(W\), and let \(b\) be the leftmost peripheral point of \(W\) with \(b\ge a+\rho_0\), if there is one. Then there is no peripheral point in \([a+\rho_0,b)\), by the choice of \(b\). Every peripheral point of \(W\) that is at least \(b\) lies in \([b,b+\rho_0)\): a point \(c\ge b+\rho_0\) would give the triple \(a<b<c\). So all peripheral points of \(W\) lie in \([a,a+\rho_0)\cup[b,b+\rho_0)\); if \(b\) does not exist they all lie in \([a,a+\rho_0)\).

\emph{Step 3: the probability for a fixed core vertex.} Fix \(\eta\in(0,p_{10}-p_{00})\) and \(\varepsilon_1\in(0,1)\), to be chosen in Step~4, and put \(\delta=\eta r_n\). Cover the intervals of Step~2 by intervals from a grid: let \(g\) range over the points \(x-\rho_1+\delta k\), \(k=0,1,\dots\), that are at most \(x+\rho_1\); there are at most \(2\rho_1/\delta+1=2p_{10}/\eta+1\) of them. For every \(a\in W\) there is a grid point \(g\) with \(g\le a<g+\delta\), so that
\[
[a,a+\rho_0)\subseteq[g,g+\rho_0+\delta).
\]
Thus if no triple exists, all peripheral points of \(W\) lie in the union of two intervals of length \(\rho_0+\delta\) with left endpoints on the grid; there are at most \(K_\eta=(2p_{10}/\eta+1)^2\) such unions. For each of them, the part of \(W\) outside it has length at least
\[
2\rho_1-2(\rho_0+\delta)=2r_n(p_{10}-p_{00}-\eta)>0 .
\]
Now condition on the labels with \(N\le\varepsilon_1n\), where \(N\) is the number of core vertices, and on \(X_i=x\). The at least \((1-\varepsilon_1)n\) peripheral positions are i.i.d.\ uniform on \([0,1]\), independent of \(x\), so all of them avoid a given set of length \(\ell\) inside \(W\) with probability at most \(e^{-(1-\varepsilon_1)n\ell}\). By a union bound over the \(K_\eta\) possible unions,
\[
\P\bigl(i\notin\widehat C_{\mathrm{star}}\mid\text{labels},X_i=x\bigr)\ \le\ K_\eta\exp\bigl(-2(1-\varepsilon_1)(p_{10}-p_{00}-\eta)nr_n\bigr).
\]

\emph{Step 4: the union bound and the choice of constants.} By Markov's inequality,
\[
\P(N>\varepsilon_1n)\le\frac{\E N}{\varepsilon_1n}=\frac{p_n}{\varepsilon_1}\to0,
\]
because \(p_n\to0\). Summing Step~3 over the vertices, each of which is a core vertex with probability \(p_n\),
\[
\P\bigl(\exists i\in C:\ X_i\in I,\ i\notin\widehat C_{\mathrm{star}}\bigr)\ \le\ \frac{p_n}{\varepsilon_1}+np_nK_\eta\exp\bigl(-2(1-\varepsilon_1)(p_{10}-p_{00}-\eta)nr_n\bigr).
\]
Assume
\[
nr_n\ge\frac{(1+\varepsilon)\log(np_n)}{2(p_{10}-p_{00})}.
\]
Since
\[
(1-\varepsilon_1)\Bigl(1-\frac{\eta}{p_{10}-p_{00}}\Bigr)\to1
\]
as \(\varepsilon_1,\eta\to0\), we can choose \(\varepsilon_1\) and \(\eta\) so small that
\[
(1-\varepsilon_1)\Bigl(1-\frac{\eta}{p_{10}-p_{00}}\Bigr)(1+\varepsilon)\ge1+\frac\varepsilon2.
\]
Then the exponent is at least \((1+\varepsilon/2)\log(np_n)\), and the second term is at most
\[
K_\eta(np_n)^{-\varepsilon/2}\to0 .
\]
Hence, with probability tending to one, every core vertex with position in \(I\) belongs to \(\widehat C_{\mathrm{star}}\). The expected number of core vertices outside \(I\) is at most \(np_n\cdot2\rho_1=2p_{10}np_nr_n\to0\), so w.h.p.\ there is none. Therefore \(C\subseteq\widehat C_{\mathrm{star}}\) w.h.p., and \(\widehat C_{\mathrm{star}}\subseteq C\) by Lemma~\ref{lem:centres}, that is, \(\P(\widehat C_{\mathrm{star}}=C)\to1\). This proves part~(i). \hfill\(\blacksquare\)

\paragraph{Part (ii).} This is Theorem~\ref{thm:sharp}(i) with \(m=1\): then \(V_1=\pi^{1/2}/\Gamma(\frac32)=2\), so
\[
c_2=V_1(p_{11}-p_{00})=2(p_{11}-p_{00}),\qquad\text{that is,}\qquad \frac1{c_2}=\frac{1}{2(p_{11}-p_{00})}.
\]

\subsection{Proof of Proposition~\ref{prop:degree}}

Write \(\rho_0=r_np_{00}\), \(\rho_1=r_np_{10}\), \(V_m=\pi^{m/2}/\Gamma(m/2+1)\) and
\[
\mu_0=nV_m\rho_0^m .
\]
The proof has two parts. First, the threshold \(\tau_n\) cannot be smaller than the typical degree of a peripheral vertex far from the boundary, since otherwise the rule would contain peripheral vertices. Second, core vertices at positions with small density \(\varphi(x)\) have a degree below this level, so the rule misses them. Both parts use that, given its position and label, the degree of a vertex is binomial (Step~D1).

\emph{Part 1: a lower bound on \(\tau_n\).} Consider a peripheral vertex \(i\) with \(X_i\in[\rho_0,1-\rho_0]^m\), an event of probability at least \(2^{-m}\) for each vertex. Its degree is Binomial\((n-1,q)\) with
\[
q\ge(1-p_n)V_m\rho_0^m .
\]
By the Chernoff bound, each such vertex has degree at least
\[
(1-\eta)(1-p_n)(n-1)V_m\rho_0^m
\]
with probability tending to one, and there is at least one such vertex with probability tending to one. If \(\tau_n\) were smaller than this value, the rule would contain a peripheral vertex with probability bounded away from zero. Hence, for large \(n\),
\[
\tau_n\ \ge\ (1-\eta)(1-p_n)(n-1)V_m\rho_0^m\ \ge\ (1-2\eta)\mu_0 .
\]

\emph{Part 2: core vertices at low-density positions.} Let \(x\) be a position with \(\varphi(x)(p_{10}/p_{00})^m\le1-5\eta\). A core vertex at \(x\) has degree Binomial\((n-1,q)\) with
\[
q\ \le\ V_m\rho_1^m\varphi(x)+p_nV_m(r_np_{11})^m\ \le\ V_m\rho_0^m\Bigl(\varphi(x)\Bigl(\frac{p_{10}}{p_{00}}\Bigr)^m+o(1)\Bigr),
\]
so its mean is at most \((1-4\eta)\mu_0\) for large \(n\). The upper-tail Chernoff bound gives
\[
\P\bigl(\deg>(1-2\eta)\mu_0\bigr)\ \le\ e^{-c\eta^2\mu_0}\ \to\ 0 ,
\]
so, with probability tending to one, the core vertex is not in the rule. \hfill\(\square\)

\section*{Acknowledgment}
The author conceptualised the research topic. The author acknowledges the interactive assistance of Anthropic's Claude Sonnet 5.5 in detailing the technical proofs and refining the language. The author accepts full accountability for the results presented herein.


\begin{thebibliography}{99}
\bibitem{Abb18} \textsc{Abbe, E.} (2018). Community detection and stochastic block models: recent developments. \emph{J. Mach. Learn. Res.} \textbf{18}(177) 1--86.
\bibitem{ABH16} \textsc{Abbe, E.}, \textsc{Bandeira, A.~S.} and \textsc{Hall, G.} (2016). Exact recovery in the stochastic block model. \emph{IEEE Trans. Inform. Theory} \textbf{62}(1) 471--487.
\bibitem{AR02} \textsc{Appel, M.~J.~B.} and \textsc{Russo, R.~P.} (2002). The connectivity of a graph on uniform points on \([0,1]^d\). \emph{Statist. Probab. Lett.} \textbf{60} 351--357.
\bibitem{ACV14} \textsc{Arias-Castro, E.} and \textsc{Verzelen, N.} (2014). Community detection in dense random networks. \emph{Ann.\ Statist.} \textbf{42} 940--969.
\bibitem{BBCH21} \textsc{Bet, G.}, \textsc{Bogerd, K.}, \textsc{Castro, R.~M.} and \textsc{van der Hofstad, R.} (2020). Detecting a botnet in a network. \emph{Math. Stat. Learn.} \textbf{3}(3/4) 315--343.
\bibitem{BE00} \textsc{Borgatti, S.~P.} and \textsc{Everett, M.~G.} (2000). Models of core/periphery structures. \emph{Social Networks} \textbf{21} 375--395.
\bibitem{BDER16} \textsc{Bubeck, S.}, \textsc{Ding, J.}, \textsc{Eldan, R.} and \textsc{R\'acz, M.~Z.} (2016). Testing for high-dimensional geometry in random graphs. \emph{Random Structures Algorithms} \textbf{49} 503--532.
\bibitem{CLWU13} \textsc{Csermely, P.}, \textsc{London, A.}, \textsc{Wu, L.-Y.} and \textsc{Uzzi, B.} (2013). Structure and dynamics of core/periphery networks. \emph{J. Complex Netw.} \textbf{1} 93--123.
\bibitem{CRLP16} \textsc{Cucuringu, M.}, \textsc{Rombach, P.}, \textsc{Lee, S.~H.} and \textsc{Porter, M.~A.} (2016). Detection of core--periphery structure in networks using spectral methods and geodesic paths. \emph{European J. Appl. Math.} \textbf{27} 846--887.
\bibitem{DC02} \textsc{Dall, J.} and \textsc{Christensen, M.} (2002). Random geometric graphs. \emph{Phys. Rev. E} \textbf{66} 016121.
\bibitem{DDC23} \textsc{Duchemin, Q.} and \textsc{De Castro, Y.} (2023). Random geometric graph: some recent developments and perspectives. In \emph{High Dimensional Probability IX} (R.~Adamczak, N.~Gozlan, K.~Lounici and M.~Madiman, eds.). \emph{Progress in Probability} \textbf{80}. Birkh\"auser, Cham. \url{https://doi.org/10.1007/978-3-031-26979-0_14}.
\bibitem{GMPS23} \textsc{Galhotra, S.}, \textsc{Mazumdar, A.}, \textsc{Pal, S.} and \textsc{Saha, B.} (2023). Community recovery in the geometric block model. \emph{J. Mach. Learn. Res.} \textbf{24}(338) 1--53.
\bibitem{Gil61} \textsc{Gilbert, E.~N.} (1961). Random plane networks. \emph{J. Soc. Indust. Appl. Math.} \textbf{9} 533--543.
\bibitem{GRK05} \textsc{Goel, A.}, \textsc{Rai, S.} and \textsc{Krishnamachari, B.} (2005). Monotone properties of random geometric graphs have sharp thresholds. \emph{Ann. Appl. Probab.} \textbf{15} 2535--2552.
\bibitem{GK00} \textsc{Gupta, P.} and \textsc{Kumar, P.~R.} (2000). The capacity of wireless networks. \emph{IEEE Trans. Inform. Theory} \textbf{46}(2) 388--404.
\bibitem{HRP08} \textsc{Higham, D.~J.}, \textsc{Ra\v{s}ajski, M.} and \textsc{Pr\v{z}ulj, N.} (2008). Fitting a geometric graph to a protein--protein interaction network. \emph{Bioinformatics} \textbf{24}(8) 1093--1099.
\bibitem{JDP83} \textsc{Joag-Dev, K.} and \textsc{Proschan, F.} (1983). Negative association of random variables with applications. \emph{Ann. Statist.} \textbf{11}(1) 286--295.
\bibitem{KM17} \textsc{Kojaku, S.} and \textsc{Masuda, N.} (2017). Finding multiple core--periphery pairs in networks. \emph{Phys. Rev. E} \textbf{96} 052313.
\bibitem{Kru96} \textsc{Krugman, P.} (1996). \emph{The Self Organizing Economy}. Wiley, Hoboken.
\bibitem{LPW09} \textsc{Levin, D.~A.}, \textsc{Peres, Y.} and \textsc{Wilmer, E.~L.} (2009). \emph{Markov Chains and Mixing Times}. American Mathematical Society, Providence.
\bibitem{LS17} \textsc{Lordan, O.} and \textsc{Sallan, J.~M.} (2017). Analyzing the multilevel structure of the European airport network. \emph{Chinese J. Aeronaut.} \textbf{30} 554--560.
\bibitem{LS19} \textsc{Lordan, O.} and \textsc{Sallan, J.~M.} (2019). Core and critical cities of global region airport networks. \emph{Phys. A} \textbf{513} 724--733.
\bibitem{MLS16} \textsc{Magone, J.~M.}, \textsc{Laffan, B.} and \textsc{Schweiger, C.} (2016). \emph{Core--Periphery Relations in the European Union: Power and Conflict in a Dualist Political Economy}. Routledge, London.
\bibitem{ML23} \textsc{Miao, R.} and \textsc{Li, X.} (2023). Informative core identification in complex networks. \emph{J. R. Stat. Soc. Ser. B} \textbf{85}(1) 108--126.
\bibitem{MNS15} \textsc{Mossel, E.}, \textsc{Neeman, J.} and \textsc{Sly, A.} (2015). Reconstruction and estimation in the planted partition model. \emph{Probab. Theory Related Fields} \textbf{162} 431--461.
\bibitem{Pen99} \textsc{Penrose, M.~D.} (1999). On \(k\)-connectivity for a geometric random graph. \emph{Random Structures Algorithms} \textbf{15} 145--164.
\bibitem{Penrose03} \textsc{Penrose, M.} (2003). \emph{Random Geometric Graphs}. Oxford University Press, Oxford.
\bibitem{RPFM14} \textsc{Rombach, M.~P.}, \textsc{Porter, M.~A.}, \textsc{Fowler, J.~H.} and \textsc{Mucha, P.~J.} (2014). Core--periphery structure in networks. \emph{SIAM J.\ Appl.\ Math.} \textbf{74} 167--190.
\bibitem{SCL20} \textsc{Sedita, S.~R.}, \textsc{Caloffi, A.} and \textsc{Lazzeretti, L.} (2020). The invisible college of cluster research: a bibliometric core--periphery analysis of the literature. \emph{Industry and Innovation} \textbf{27} 562--584.
\bibitem{SSH18} \textsc{Solovey, K.}, \textsc{Salzman, O.} and \textsc{Halperin, D.} (2018). New perspective on sampling-based motion planning via random geometric graphs. \emph{Int. J. Robot. Res.} \textbf{37}(10) 1117--1133.
\bibitem{Tsy09} \textsc{Tsybakov, A.~B.} (2009). \emph{Introduction to Nonparametric Estimation}. Springer, New York.
\bibitem{WPKWC22} \textsc{Wedell, E.}, \textsc{Park, M.}, \textsc{Korobskiy, D.}, \textsc{Warnow, T.} and \textsc{Chacko, G.} (2022). Center--periphery structure in research communities. \emph{Quant. Sci. Stud.} \textbf{3} 289--314.
\bibitem{YS23} \textsc{Yanchenko, E.} and \textsc{Sengupta, S.} (2023). Core--periphery structure in networks: a statistical exposition. \emph{Stat. Surv.} \textbf{17} 42--74.
\bibitem{YS26} \textsc{Yanchenko, E.} and \textsc{Sengupta, S.} (2026). A label-switching algorithm for fast core--periphery identification. \emph{Netw. Sci.} \textbf{14}.
\bibitem{YSM25} \textsc{Yanchenko, E.}, \textsc{Sengupta, S.} and \textsc{Mukherjee, S.} (2025). Statistical inference for core--periphery structures. Preprint, arXiv:2508.04730.
\bibitem{Yuan24} \textsc{Yuan, M.} (2024). Asymptotic distribution of the friendship paradox of a random geometric graph. \emph{Braz. J. Probab. Stat.} \textbf{38}(3). \url{https://doi.org/10.1214/24-BJPS615}.
\bibitem{Yuan26} \textsc{Yuan, M.} (2026). Asymptotic distribution of the global clustering coefficient in a random annulus graph. \emph{Acta Math. Hungar.} \url{https://doi.org/10.1007/s10474-026-01644-5}.
\bibitem{YuanWLLN26} \textsc{Yuan, M.} (2026). The weak law of large numbers for the friendship paradox index. \emph{ALEA Lat. Am. J. Probab. Math. Stat.}, accepted. \url{https://arxiv.org/abs/2602.10055}.
\bibitem{YY26} \textsc{Yuan, M.} and \textsc{Yu, F.} (2026). Hypothesis testing for the dimension of random geometric graphs. \emph{Electron. J. Stat.} \textbf{20}(2) 5431--5469. \url{https://doi.org/10.1214/26-EJS2609}.
\bibitem{ZMN15} \textsc{Zhang, X.}, \textsc{Martin, T.} and \textsc{Newman, M.~E.~J.} (2015). Identification of core--periphery structure in networks. \emph{Phys. Rev. E} \textbf{91} 032803.
\bibitem{Zha15} \textsc{Zhao, J.} (2015). The absence of isolated node in geometric random graphs. In \emph{53rd Annual Allerton Conference on Communication, Control, and Computing (Allerton)} 881--886. Monticello, IL.
\end{thebibliography}
\end{document}